\documentclass{article}
\usepackage{graphicx} % Required for inserting images
\usepackage{geometry}               % See geometry.pdf to learn the layout options. There are lots.
\usepackage[usenames,dvipsnames]{xcolor}
\usepackage{amssymb,amsmath,amsthm,amsfonts}
\usepackage{bm}                     % mathbf greek
\usepackage[english]{babel}
\usepackage[utf8]{inputenc}
\usepackage[colorlinks]{hyperref}
\usepackage{comment}
\usepackage[showonlyrefs]{mathtools}
\mathtoolsset{showonlyrefs=true}
\usepackage{enumitem}
\newtheorem{proposition}{Proposition}[section]
\newtheorem{theorem}[proposition]{Theorem}

\newtheorem{lemma}[proposition]{Lemma}

\theoremstyle{definition}
\newtheorem{definition}[proposition]{Definition}
\theoremstyle{definition}
\newtheorem{hypothesis}{Hypothesis}

\theoremstyle{remark}
\newtheorem{remark}[proposition]{Remark}
\numberwithin{equation}{section}
\theoremstyle{definition}

\newcommand{\abs}[1]{\left|#1\right|}
\newcommand{\norm}[1]{\left\lVert#1\right\rVert}

\newcommand{\eps}{\varepsilon}

\newcommand{\N}{{\mathbb{N}}}

\newcommand{\R}{{\mathbb{R}}}

\DeclareMathOperator*{\argmin}{argmin}

\DeclareMathOperator{\dive}{div}
\DeclareMathOperator{\ent}{Ent}

\newcommand{\AC}[1]{{\color{violet} [AC: {#1}]}}

\usepackage[maxbibnames=9, giveninits=true, doi=false, url=false, isbn=false]{biblatex}
\usepackage{doi}
\usepackage{mathrsfs}
\title{Small noise limit for Entropic Variational Mean-Field Games}
\author{Alessandro Cosenza\thanks{Universit\'e Paris-Saclay, CNRS, Laboratoire de math\'ematiques d'Orsay, ParMA, Inria Saclay, 91405, Orsay, France. Email: \href{mailto:alessandro.cosenza@universite-paris-saclay.fr}{\texttt{alessandro.cosenza@universite-paris-saclay.fr}}}
\and Luca Nenna\thanks{Universit\'e Paris-Saclay, CNRS, Laboratoire de math\'ematiques d'Orsay, ParMA, Inria Saclay, 91405, Orsay, France. Email: \href{mailto:luca.nenna@universite-paris-saclay.fr}{\texttt{luca.nenna@universite-paris-saclay.fr}}}\ \thanks{Institut Universitaire de France (IUF).}}
\begin{document}

\maketitle
\begin{abstract}
    We study second order variational mean-field games with quadratic Hamiltonian. In analogy with entropic optimal transport, Benamou, Carlier, Di Marino and the second author recently proposed an entropy-minimization formulation of the problem, which could be interpreted as the Lagrangian counterpart of the classical Eulerian approach. We investigate the small-noise $\Gamma$-limit for the Lagrangian entropy-minimization formulation and establish $\Gamma$-convergence of second order mean-field games towards their deterministic first order counterpart. This also translates into another (dynamical) proof of the convergence of the Schr\"odinger problem to optimal transport. Moreover, we study the original Eulerian problem and prove a first order energy expansion, which establishes a selection principle for minimizers of the deterministic problem.
\end{abstract}
\noindent
{\footnotesize \textbf{AMS-Subject Classification}}. 49N80, 49J45, 49Q22, 91A16.
{\footnotesize }
    \tableofcontents
\section{Introduction}

Mean-Field Games were independently introduced in the seminal papers by Lasry and Lions \cite{MFGLL1,MFGLL2,MFGLL3} and Huang, Caines and Malham\'e \cite{huang_malhame_caines_2006}, with the aim of describing Nash equilibria in differential games with infinitely many indistinguishable agents. The evolution of the agents is formally described by a coupled Hamilton--Jacobi and Fokker--Planck system. In the setting of a decoupled Hamiltonian with quadratic nonlinearity, this system takes the form
\begin{equation}\label{e:thesystem}
\left\{\begin{aligned}
-\partial_t u-\frac{\nu}{2}\Delta u
  +\frac12|\nabla u|^2 &= f[\rho_t],
  &&\text{in }(0,1)\times\mathbb R^d,\\
\partial_t\rho-\frac{\nu}{2}\Delta\rho
  -\operatorname{div}(\rho\nabla u)&=0,
  &&\text{in }(0,1)\times\mathbb R^d,\\
\rho_{|t=0}=\rho_0,\qquad
u_{|t=1}&=g[\rho_1].
\end{aligned}\right.
\end{equation}
Here $\rho$ represents the distribution of the agents, while $u$ is a value function representing the strategy of the agents. The first equation describes optimal individual behaviour given the evolution of the population, while the second equation transports the population with a drift term which depends on the optimal strategy. The number $\nu\geq 0$ represents the diffusivity of the population. When $\nu>0$ the population is subject to diffusion and the problem is referred to as second-order Mean-Field Games, while for $\nu=0$ one talks about first-order Mean-Field Games. 
The terms $f$ and $g$ represent the coupling of the system.

We focus on \emph{potential} Mean-Field Games, for which the couplings $f$ and $g$ arise, formally, as first variations of interaction energies $F$ and $G$. In this case, the PDE system is associated with a variational problem \cite{MFGLL1,MFGLL2,MFGLL3,variationalMFG}. If $v$ denotes the population velocity, solutions to \eqref{e:thesystem} can be found as minimizers of the functional
\begin{equation}
\label{e:introeul}
\mathcal J(\rho,v)
=\frac12\int_0^1\!\int_{\mathbb R^d}|v_t|^2\,d\rho_t\,dt
 +\int_0^1 F(\rho_t)\,dt+G(\rho_1)
\end{equation}
under the constraint
\begin{equation}\label{e:intro-fp}
\partial_t\rho-\frac{\nu}{2}\Delta\rho
  +\operatorname{div}(\rho v)=0,
  \qquad \rho_{|t=0}=\rho_0.
\end{equation}
When $\nu=0$, problem \eqref{e:introeul}-\eqref{e:intro-fp} resembles the Benamou--Brenier formula \cite{Benamoubrenier} for the optimal transportation problem. The main difference of \eqref{e:introeul} is that only one marginal is fixed and some running and terminal costs, represented by $F$ and $G$, are present.  Optimal transport can be formally recovered by choosing $F=0$ and $G=+\infty\chi_{\{\rho_1\}^c}$, see Section \ref{sec:hyp} for the precise set of hypotheses we use.
Analogously, for $\nu>0$, a parallel can be found with the Schr\"odinger bridge problem, an entropy minimization problem which can be interpreted as a regularization of optimal transport, see \cite{survey,leonard_convergence} and \cite{Nutzintr} for an introduction to the topic. In particular, problem \eqref{e:introeul}-\eqref{e:intro-fp} is similar to the ``noisy" Benamou--Brenier formula for entropic optimal transport, which has received attention in recent years (see \cite{leonardbb,noisyBB,giglitam} and also the recent paper \cite{mattia} for a PDE approach to the formula).

Both the standard and the ``noisy" Benamou--Brenier formulas provide an \emph{Eulerian} dynamical reformulation of the two static problems of classical and entropic optimal transport. Another dynamical approach to such problems is the \emph{Lagrangian} one, that is considering  minimization problems on the path space $\mathcal{P}(C([0,1],\R^d))$. For classical optimal transport, the kinetic energy of the paths is minimized (see any classical reference for OT, for instance \cite{Villani,AGS}), while for entropic optimal transport, the relative entropy with respect to a reversible Wiener measure is minimized \cite{survey}. Going back to the Mean--Field Game problem, the presence of a running cost makes the problem purely dynamic, so that no static formulation is available. However, one may look for a Lagrangian reformulation of  \eqref{e:introeul}-\eqref{e:intro-fp}. When $\nu=0$, such a model has been proposed in \cite{variationalMFG} and the equivalence is essentially due to classical optimal transport results. When $\nu>0$ the situation is more complicated, and the corresponding result was proved by
 Benamou, Carlier, Di Marino and the second author in \cite{benamou2019entropy}. Building on some ideas of \cite{dawsongartner}, they proposed an entropy-minimization problem on path space equivalent to its Eulerian counterpart, and a corresponding Sinkhorn algorithm.

 The aim of this paper is to investigate the small-noise behavior of both the Eulerian and Lagrangian descriptions, in the sense of $\Gamma$-convergence. Our results are in the same spirit as the small noise $\Gamma$-convergence results of entropic optimal transport such as \cite{leonard_convergence,carlier2017convergence}. Those papers, however, leverage the static formulation available for optimal transport, while for our problem we have to rely on a purely dynamical approach. To the best of our knowledge, the only paper which proposes a dynamical approach to the small noise $\Gamma$-limit of entropic transport is \cite{baradatm}, whose techniques we adapt to our setting.  
As far as the convergence of minimizers is concerned, our $\Gamma$-convergence results are qualitative, in contrast to recent works on convergence rates for vanishing-viscosity second-order Mean-Field Games \cite{rate1,rate2}. On the other hand, the $\Gamma$-convergence approach only requires weak regularity and structural assumptions on the data, and it directly yields the convergence of minimum values and of minimizers.

Let us point out the main difficulties. First, the presence of a running cost makes the problem purely dynamical, so that the static arguments available for entropic optimal transport \cite{leonard_convergence,carlier2017convergence} cannot be used: recovery sequences have to be built directly at the level of curves of measures and of measures on paths. Second, we work on the whole space $\R^d$ and the reference measure $R^\nu$ is only $\sigma$-finite, so that compactness is not automatic; in the Lagrangian setting it is recovered from the exponential tightness of Brownian motion. Third, the natural Lagrangian recovery sequence, obtained by perturbing each path with a Brownian bridge, leaves the endpoint coupling unchanged, so that its entropic cost does not vanish in the limit: the endpoint coupling has to be regularized while keeping control of the kinetic energy. Finally, in the Mean-Field Game setting the terminal marginal is free and may have infinite entropy, which requires an additional ad hoc regularization, and which also prevents us from going beyond a first-order energy expansion.

To detail our results, let us introduce the models at play; see Section \ref{sec:hyp} for the precise definitions. 
 For a curve $\rho\in C([0,1];\mathcal P_2(\mathbb R^d))$, let
\[
\mathcal{A}(\rho)=\inf\left\{
 \frac12\int_0^1\!\int_{\mathbb R^d}|v_t|^2\,d\rho_t\,dt:
 \partial_t\rho+\operatorname{div}(\rho v)=0
\right\}.
\]
 Up to a change of variable (see \cite{benamou2019entropy}) problem \eqref{e:introeul}-\eqref{e:intro-fp} is equivalent to minimizing the functional 
\begin{equation}
\label{e:intro-Eulerian-functional}
\mathcal F_\nu(\rho)
=\mathcal{A}(\rho)+\int_0^1 F(\rho_t)\,dt+G(\rho_1)
 +\frac{\nu}{2}\bigl(\operatorname{Ent}(\rho_1)-\operatorname{Ent}(\rho_0)\bigr)
 +\frac{\nu^2}{8}\int_0^1 I(\rho_t)\,dt,
\end{equation}
subject to $\rho_{|t=0}=\rho_0$. Here $\operatorname{Ent}(\rho)=\int\rho\log\rho\,dx$ denotes the entropy (i.e.\ minus the Shannon entropy) and $I$ the Fisher information. 

For the Lagrangian description, let $\Omega=C([0,1];\mathbb R^d)$ and $Q\in\mathcal P(\Omega)$. Its marginals are $Q_t=(e_t)_\#Q$, with $e_t(\omega)=\omega_t$. Denote by $R^\nu$ the reversible Wiener measure with diffusivity $\nu$ (see Section \ref{sec:wiener}). For $\nu>0$, the entropy-minimization approach of \cite{benamou2019entropy} is
\begin{equation}\label{e:intro-Lagrangian-functional}
\mathcal L_\nu(Q)
 =\nu H(Q\mid R^\nu)
 +\int_0^1 F(Q_t)\,dt+G(Q_1),
 \qquad Q_0=\rho_0,
\end{equation}
where $H(Q\mid R^\nu)$ represents relative entropy of $Q$ with respect to $R^\nu$ (a $\sigma$-finite measure), see Section \ref{sec:entr}. The corresponding limit functional is
\begin{equation}\label{e:intro-Lagrangian-limit}
\mathcal L_0(Q)=\mathcal{K}(Q)+\int_0^1F(Q_t)\,dt+G(Q_1),
\qquad
\mathcal{K}(Q)=\int_\Omega k(\omega)\,dQ(\omega),
\end{equation}
with $Q_0=\rho_0$, where
\[
k(\omega)=
\begin{cases}
\displaystyle \frac12\int_0^1|\dot\omega_t|^2\,dt,
&\omega\in AC^2([0,1];\mathbb R^d),\\
+\infty,&\text{otherwise.}
\end{cases}
\]
For $\nu>0$, the results of \cite{benamou2019entropy} establish that
the Eulerian and Lagrangian minimization problems \eqref{e:intro-Eulerian-functional} and \eqref{e:intro-Lagrangian-functional} are equivalent up to an additive constant, while the equivalence of \eqref{e:intro-Eulerian-functional} and \eqref{e:intro-Lagrangian-limit} for $\nu=0$ is established in \cite{variationalMFG}.
We make the following assumptions on the data $\rho_0,F,G$ (see Section \ref{sec:hyp} for details):
\begin{itemize}
\item [(\ref{hyp:initial})] $\rho_0\in\mathcal P_2(\mathbb R^d)$ has finite entropy;

\item [(\ref{hyp:mfg})] For  the   Mean--Field Games framework  described above: $F$ and $G$ are the sum of convex local energies with polynomial growth and nonlocal $W_2$ continuous terms (in the language of Mean--Field Games, we suppose that both \emph{local} and \emph{non-local} couplings are present);

\item[(\ref{hyp:eot})] For the entropic optimal transport framework, corresponding to a fixed terminal marginal with no running cost: $F=0$ and $G$ is the indicator of a prescribed $\rho_1\in\mathcal P_2(\mathbb R^d)$ which we assume has finite entropy.
\end{itemize}  

We now state our main results. Our first result is the Eulerian $\Gamma$-limit (see Theorem \ref{th:maineul}).

\begin{theorem}[Eulerian $\Gamma$-convergence]\label{th:intro-eulerian}
Suppose \ref{hyp:initial} and either \ref{hyp:mfg} or \ref{hyp:eot} hold. 
As $\nu\downarrow0$, the functionals $\mathcal F_\nu$ $\Gamma$-converge to $\mathcal F_0$ with respect to uniform narrow convergence in $C^{0}([0,1],\mathcal{P}(\R^d))$ (see Section \ref{sec:pre}).
\end{theorem}
Secondly, we establish a first order energy expansion in the spirit of \cite{giglitam,confortitam,legerexpansion}, in the setting of Mean--Field Games. Set
\[
\mathcal{ F}^1(\rho)=\frac12\bigl(\operatorname{Ent}(\rho_1)
 -\operatorname{Ent}(\rho_0)\bigr),\qquad
\mathcal{F}^2(\rho)=\frac18\int_0^1I(\rho_t)\,dt,
\]
so that $\mathcal F_\nu=\mathcal F_0+\nu\mathcal {F}^1+\nu^2\mathcal{F}^2$. Denote
\[
m_\nu=\inf\mathcal F_\nu,\qquad
\mathcal U_0=\operatorname{argmin}\mathcal F_0,
\qquad
m_1=\inf_{\rho\in\mathcal U_0}\mathcal{F}^1(\rho),
\qquad
\mathcal U_1=\operatorname{argmin}_{\rho\in\mathcal U_0}\mathcal{F}^1(\rho).
\]
We prove the following energy expansion and selection principle (see Theorem \ref{th:energyexp}).

\begin{theorem}[First-order expansion and selection principle]\label{thm:intro-selection}
Suppose \ref{hyp:initial} and \ref{hyp:mfg} hold. Suppose that $m_1<+\infty$ and that there exists $\bar\rho\in\mathcal U_1$ for which $\mathcal{F}^2(\bar\rho)<+\infty$. Then
\begin{equation}\label{e:intro-minimum-expansion}
\boxed{\,m_\nu=m_0+\nu m_1+o(\nu)\quad\text{as }\nu\rightarrow 0\,}
\end{equation}
Moreover, if $\nu_n\rightarrow 0$ and $\rho^n$ is a minimizer of $\mathcal F_{\nu_n}$, then every cluster point of $(\rho^n)_n$ belongs to $\mathcal U_1$.
\end{theorem}

Finally we prove the Lagrangian $\Gamma$-limit (see Theorem \ref{th:mainlagr}). This is the main result of the paper. 

\begin{theorem}[Lagrangian $\Gamma$-convergence]\label{thm:intro-lagrangian}
Suppose \ref{hyp:initial} and either \ref{hyp:mfg} or \ref{hyp:eot} hold. As $\nu\rightarrow 0$, the functionals $\mathcal L_\nu$ $\Gamma$-converge to $\mathcal L_0$ with respect to narrow convergence on $\mathcal{P}(\Omega)$.
\end{theorem}
Combined with the equicompactness results of Propositions \ref{prop:compeul} and \ref{prop:complagr}, Theorems \ref{th:intro-eulerian} and \ref{thm:intro-lagrangian} yield the convergence of the minimum values and, up to subsequences, of the minimizers. In particular, since $\mathcal L_\nu$ is the functional addressed numerically by the Sinkhorn-type algorithm of \cite{benamou2019entropy}, Theorem \ref{thm:intro-lagrangian} shows that this entropic regularization provides a consistent approximation of first-order Mean-Field Game equilibria, at the level of the distribution of trajectories.

We now comment on each theorem and on the strategy of the proofs.
The validity of Theorem \ref{th:intro-eulerian} is expected and essentially folklore, but we could not locate a proof in the literature. Its proof is based on \cite[Theorem 2.3]{baradatm}, where they prove the analogous result for entropic optimal transport on the torus. We adapt their strategy to the setting of \ref{hyp:mfg}, with the added difficulty of working on the whole space $\R^d$. Since we work on an unbounded domain, we first establish, as usual for $\Gamma$-convergence results, an equicompactness result (see Proposition \ref{prop:compeul} and Remark \ref{rem:comp}) which motivates the choice of the topology for the $\Gamma$-limit. The proof of the $\Gamma$-$\liminf$
inequality is easy once the problem is formulated as in \eqref{e:intro-Eulerian-functional}, since the sequences of curves satisfy the same constraint as the limit problem (the continuity equation), and the interaction terms are lower-semicontinuous. For the $\Gamma$-$\limsup$ inequality, the main ingredient is Lemma \ref{lem:limsup}, in which we adapt the computations of \cite[Lemma 4.1]{baradatm} to prove that the recovery sequence given by an appropriate convolution with the heat semigroup $\tau_\nu$, that is
\[
\rho_t^\nu=\rho_t*\tau_{\nu t},
\]
has kinetic energy controlled by the original curve. Since the convergence of the recovery sequence is stronger (namely in $C^0([0,1],\mathcal{P}_2(\R^d))$), this is enough to prove convergence of the interaction terms.

Let us now get to Theorem \ref{thm:intro-selection}. Our goal here is to establish a result similar in spirit to second order energy expansions for entropic optimal transport \cite{giglitam,confortitam,legerexpansion}. Since we allow for non-convex (e.g.\ nonlocal) interaction terms, uniqueness of minimizers of the limit problem is not expected (see \cite{bardi} for an overview of non-uniqueness in Mean--Field Games). Hence, proving the energy expansion \eqref{e:intro-minimum-expansion} also establishes a selection principle for converging sequences of minimizers (see also \cite{louis}, where a similar question is studied in the case of non-uniqueness of transport plans).  When the limit problem has more than one minimizer, the ones with the smallest value of $\operatorname{Ent}(\rho_1)$ are selected, see Section \ref{sec:expansion} for an example of application. In this example, the population starts from a Gaussian and is attracted at the final time either towards a translated or towards a dilated Gaussian: the first-order problem has two minimizers with the same cost, and the noise selects the one which spreads the population, i.e.\ the one with smaller terminal entropy. Even though we do not state it explicitly, this could be interpreted in terms of asymptotic developments by $\Gamma$-convergence as in \cite{higherordergammaconv}. The proof of Theorem \ref{thm:intro-selection} follows similar lines to \cite{legerexpansion,confortitam}. We just remark that, in contrast with the entropic optimal transport, we only get a first order development. Indeed, unlike in optimal transport, in which both marginals are fixed, the terminal value of a Mean-Field Game is generally free to vary, which makes a second-order argument not immediate (see Remark \ref{rem:exp}).

Finally we comment on Theorem \ref{thm:intro-lagrangian}. We first establish the existence of a minimizer for the Lagrangian problem under our set of hypotheses (Proposition \ref{prop:existencelagr}), which differs slightly from the one in \cite{benamou2019entropy}. Then, as in the Eulerian case, we prove an equicompactness statement (Proposition \ref{prop:complagr}). This is not immediate since we work in an unbounded domain and the reversible Wiener measure $R^\nu$ is only $\sigma$-finite. Compactness is recovered by exploiting the exponential tightness of Brownian motion, a consequence of Schilder's Theorem (Lemma \ref{lem:schilder}). The proof of the $\Gamma$-limit is instead inspired by the results of \cite[Theorem 3.1]{baradatm} for the so-called Br\"odinger problem on the torus \cite{brodinger}. The $\Gamma$-$\liminf$ inequality can be directly adapted, and it is based on the duality formula for relative entropy and an approximation of the kinetic energy $\mathcal{K}$ (see Lemma \ref{lem:baradatkin}).
The $\Gamma$-$\limsup$ inequality is more involved, and is the main contribution of this paper. Here, the situation differs if we consider the sets of hypotheses \ref{hyp:mfg} or \ref{hyp:eot}. In both cases, the main idea is to construct the recovery sequence by perturbing with a Brownian bridge $R^{\nu,0,0}$ each curve which is charged by a limit path $Q\in\mathcal{P}(\Omega)$:
\begin{equation}
    Q^\nu=\int B^\nu_\omega \,dQ(\omega), \qquad B^{\nu}_\omega=(\omega+\cdot)\#R^{\nu,0,0}.
\end{equation}
This leaves the endpoint coupling $\gamma=(e_0,e_1)\#Q$ unchanged, while resulting in an estimate on the entropy (Proposition \ref{prop:entrbridge}): 
    \begin{equation}
\nu H(Q^\nu|R^{\nu})\leq \nu H(\gamma|R^{\nu}_{0,1})+\mathcal{K}(Q)-\frac{1}{2}\int_{\R^d\times\R^d}\abs{x-y}^2d\gamma,
    \end{equation}
    where $R^{\nu}_{0,1}$ is the heat transition kernel.
To conclude, we would like the first and last terms on the right hand side to vanish in the limit, but this is not possible, since the bridge construction leaves the endpoints unchanged. Hence, we perform first an appropriate regularization of the endpoint coupling $\gamma^{\nu}$, paired with a stability result of the kinetic energy $\mathcal{K}$. In the case of \ref{hyp:eot}, we choose $\gamma^{\nu}$ to be a recovery sequence for the static $\Gamma$-limit \cite{carlier2017convergence}, while the role of the stability result is played by Proposition \ref{prop:stability}, which is based on \cite{baradatcont}. This establishes (another) proof of the convergence of entropic optimal transport on $\R^d$ in the dynamical (Lagrangian) setting. In the case of \ref{hyp:mfg}, a further regularization is necessary, since in general the final marginal $Q_1$ could have infinite entropy. We thus perform an ad-hoc regularization of the endpoint coupling which also bounds the entropy of $Q_1$ (Lemma \ref{lem:moll}), and we pair it with a different stability result for $\mathcal{K}$ (Proposition \ref{prop:stability2}). Since, at the level of marginal law, the above construction amounts to convolutions with heat kernels, this is also enough to have convergence of the interaction terms.  

The rest of the paper is organized as follows. Section \ref{sec:pre} collects preliminaries on optimal transport, Wiener measures and relative entropy, and sets out the Eulerian and Lagrangian problems. Section \ref{sec:eul} proves the Eulerian $\Gamma$-limit and the first-order energy expansion. Finally, Section \ref{sec:lagr} establishes the Lagrangian $\Gamma$-limit.

\section{Preliminaries}
\label{sec:pre}
\subsection{Notation}
We write $a\lesssim_\alpha b$
if $a\leq C_\alpha b$ for a constant $C_\alpha>0$ depending only on the parameter $\alpha$.
We denote by $d\in\mathbb{N}$  the dimension of the
ambient Euclidean space. If $X$ is a Polish space, we denote by
$\mathcal{P}(X)$ the space of Borel probability measures on $X$. 
 If $T:X\to Y$ is a measurable
map and $\mu\in\mathcal{P}(X)$, then $T_\#\mu$ denotes the push-forward
of $\mu$ through $T$. For $\mu_0,\mu_1\in \mathcal{P}(X)$, we denote by $\Pi(\mu_0,\mu_1)$ the set of their couplings. We write $\mu_0\ll\mu_1$ for absolute continuity and
$\mu_0*\mu_1$ for convolution.
Set $\Omega=C([0,1];\mathbb{R}^d)$.
For $t\in[0,1]$, the evaluation map is 
    $e_t:\Omega\to\mathbb{R}^d,
    e_t(\omega)=\omega_t.$
For $Q\in\mathcal{P}(\Omega)$, its time marginal is denoted by $Q_t=(e_t)_\# Q.$
We denote by $\mathcal{N}(m,C)$ the law of a gaussian variable on $\R^d$ with mean $m$ and covariance matrix $C$.
For $s>0$, the heat kernel associated with $\frac12\Delta$ is
\[
    \tau_s(x)
    =
    \frac{1}{(2\pi s)^{d/2}}
    \exp\!\left(-\frac{|x|^2}{2s}\right),
    \qquad x\in\mathbb{R}^d.
\]
That is,  $\tau_s=\mathcal{N}(0,sI_d)$, the law of $W_s$ for a standard Brownian motion $W$. For a random variable $X$ taking values in a Polish space $E$,
we write $\operatorname{Law}(X)$ for its law, and $X\sim\mu$ if
$\operatorname{Law}(X)=\mu$.
\subsection{Optimal Transport and metrics}
Let $(X,d)$ be a Polish metric space. We denote by
\[
\mathcal P_2(X)
=
\left\{
\mu\in\mathcal P(X):
\int_X d(x,x_0)^2\,d\mu(x)<+\infty
\right\},
\]
where the definition is independent of the choice of $x_0\in X$.
For $\mu_0,\mu_1\in\mathcal P_2(X)$, we denote by $\Pi(\mu_0,\mu_1)$ the set of couplings of $\mu_0$ and $\mu_1$ and define their quadratic Wasserstein distance 
\[
    W_2^2(\mu_0,\mu_1)
 =
    \inf_{\gamma\in\Pi(\mu_0,\mu_1)}
    \int_{X\times X} d_X(x,y)^2\,d\gamma(x,y).
\]
We also define for $\mu_0,\mu_1\in\mathcal{P}(X)$, the ``cut" Wasserstein distance 
\begin{equation}   
\label{e:narrow}
W_N(\mu_0,\mu_1)=\inf_{{\gamma\in\Pi(\mu_0,\mu_1)}}\int_{X\times X}\min\{1,d_X(x,y)\}\,d\gamma.
\end{equation}
Recall that \eqref{e:narrow} metrizes narrow convergence on $\mathcal{P}(X)$, see \cite[Corollary 6.13]{Villani}.

In particular for $X=\mathbb R^d$ we let $d_{\R^d}$ be the Euclidean distance, while for $X=\Omega=C^0([0,1];\mathbb R^d)$,
we consider the distance
\[
    d_{\Omega}(\omega,\eta)
    =
    \|\omega-\eta\|_\infty
    =
    \sup_{t\in[0,1]}|\omega_t-\eta_t|.
\]
The quadratic Wasserstein distance admits the equivalent
probabilistic representation
\begin{equation}
\label{e:wassprob}
        W_2^2(\mu_0,\mu_1)
    =
    \inf\left\{
        \mathbb E[d(X,Y)^2]:
        X\sim\mu_0,\ Y\sim\mu_1
    \right\},
\end{equation}
where the infimum runs over all pairs of random variables
$(X,Y)$ defined on a common probability space, and is attained.

We endow $C^0([0,1];\mathcal P_2(\mathbb R^d))$ with 
\begin{equation}
\label{e:uniform}
    d_{C^0}(\rho,\widetilde\rho)
    =
    \sup_{t\in[0,1]}W_2(\rho_t,\widetilde\rho_t).
\end{equation}
and say that $\rho_n\rightarrow \rho$ in $C^0([0,1];\mathcal P_2(\mathbb R^d))$ if $d_{C^0}(\rho_n,\rho)\rightarrow 0$.
We also need a notion of uniform narrow convergence on $C^{0}([0,1];\mathcal{P}(\R^d))$. We say that $\rho^n\rightarrow \rho$ narrowly in $C^{0}([0,1];\mathcal{P}(\R^d))$ if 
\begin{equation}
d_N(\rho^n,\rho)=\sup_{t\in[0,1]}W_N(\rho_t^n,\rho_t)\rightarrow 0.
\end{equation}
 \begin{lemma}
  \label{lem:cont2}
    Let $Q^n,Q\in \mathcal{P}_2(\Omega)$ such that $W_2(Q^n,Q)\rightarrow 0$. Then $(Q^{n}_t)_{t\in[0,1]}\rightarrow (Q_t)_{t\in[0,1]}$ in $ C([0,1];\mathcal{P}_2(\R^d))$. If instead $Q^n\rightarrow Q$ narrowly, then $(Q^{n}_t)_{t\in[0,1]}\rightarrow (Q_t)_{t\in[0,1]}$ narrowly in $ C([0,1];\mathcal{P}(\R^d))$.
\end{lemma}
\begin{proof}
    Take $\Gamma_n\in \mathcal{P}(\Omega\times\Omega)$ optimal for $W_2(Q^n,Q)$. Define for $n\in\N$ and $t\in [0,1]$ the plan $\pi_n^t=(e_t\otimes e_t)\#\Gamma_n\in \Pi(Q^n_t,Q_t)$. For any $t\in [0,1]$, it holds
    \begin{equation}
    W_2^2(Q^n_t,Q_t)\leq\int_{\R^d\times\R^d}\abs{x-y}^2\,d\pi_n^t(x,y)= \int\abs{\gamma_t-\eta_t}^2\,d\Gamma_n(\gamma,\eta)\leq W_2^2(Q^n,Q).
    \end{equation}
    Taking the supremum in time and sending $n$ to infinity yields the first statement.
    The proof of the second statement is identical, choosing now $\Gamma_n\in \mathcal{P}(\Omega\times\Omega)$ optimal for \eqref{e:narrow}.
\end{proof}

We finally recall  the Benamou-Brenier formula \cite{Benamoubrenier}: for $\mu_0,\mu_1\in \mathcal{P}_2(\R^d)$
\begin{equation}
\label{e:bb}
    W^2_2(\mu_0,\mu_1)= \min_{(\rho,v)}\left\{\int_0^1\int_{\R^d}\abs{v}^2d\rho_t d t \colon \ \partial_t\rho +\dive (\rho v)=0, \ \rho_0=\mu_0, \ \rho_1=\mu_1 \right\},
\end{equation}
where $\rho\in C^{0}([0,1]; \mathcal{P}_2(\R^d))$, $v$ is a measurable vector field and the continuity equation holds distributionally. 
\subsection{Reversible Wiener measure}
\label{sec:wiener}
For $\nu>0$ and $x\in\mathbb R^d$,
$R^{\nu,x}\in\mathcal P(\Omega)$ denotes the law of
$(x+\sqrt\nu\,W_t)_{t\in[0,1]}$, where $W$ is a standard
$d$-dimensional Brownian motion. In particular,
\[
  (e_t)_\#R^{\nu,x}=\tau_{\nu t}(\cdot-x).
\]
The reversible Wiener measure with Lebesgue initial marginal is
\begin{equation}
\label{e:wiener}
  R^\nu=\int_{\mathbb R^d}R^{\nu,x}\,dx.
\end{equation}
Notice that $R^\nu$ is an unbounded $\sigma$-finite measure. Moreover, $R^\nu_t=dx$ for all $t\in[0,1]$.

We denote the endpoint reference measure by
\begin{equation}
  R_{0,1}^\nu=(e_0,e_1)_\#R^\nu,
  \qquad R_{0,1}^\nu(dx,dy)=\tau_\nu(y-x)\,dx\,dy.
\end{equation}
For $x,y\in\mathbb R^d$, the Brownian bridge
$R^{\nu,x,y}$ is the law of Brownian motion with diffusivity
$\nu$ conditioned to start at $x$ and end at $y$. That is, it is the measure given by the disintegration theorem (see for instance \cite[Appendix B]{figalli_glaudo}) such that
\[
  R^\nu=\int_{\mathbb R^d\times\mathbb R^d}
     R^{\nu,x,y}\,R_{0,1}^\nu(dx,dy).
\]
We write $B^\nu=R^{\nu,0,0}$. We refer for instance to \cite{revuzyor} for more properties of brownian motion and brownian bridges. We just remark that 
\begin{equation}
\label{e:bridge}
    B^\nu_t=\tau_{\nu t(1-t)}dx
\end{equation}
and
\begin{equation}
\label{e:linfest}
  \mathbb E[\sup_t\abs{B_t^\nu}^2]\le C_d\nu.   
\end{equation}
\subsection{Relative Entropy}
\label{sec:entr}
Let $X$ be a Polish space.
For $q,r \in \mathcal{P}(X)$ we define the relative entropy 
\begin{equation}
\label{e:entr}
    H(q|r)=
    \begin{cases}
    \int_{X} \log(\frac{dq}{dr})dq & \text{ if } q\ll r, \\
    +\infty & \text{ otherwise }.
    \end{cases}
\end{equation}
An equivalent representation for the relative entropy is the following
duality formula (see for instance \cite{Nutzintr}):
\begin{equation}
\label{e:duality}
    H(q|r)=\sup\left\{\int_{X}fdq-\log\int_Xe^{f}dr; f:X\rightarrow \R \text{ bounded and measurable}\right\}.
\end{equation}
We also recall the behavior of entropy with respect to disintegration (see for instance \cite[Theorem 2.4]{leonard_path}). 
\begin{proposition}
\label{prop:entrdisint}
Let $X,Y$ be Polish spaces, let $p,r\in\mathcal P(X)$.
Given a measurable map $\Phi:X\to Y$, denote by
$p^y$ and $r^y$ the corresponding
probability measures given by the disintegration theorem, so that
\[
    p=\int_Y p^y\,d(\Phi_\#p)(y),
    \qquad
    r=\int_Y r^y\,d(\Phi_\#r)(y).
\]
Then
\[
    H(p\mid r)
    =
    H(\Phi_\#p\mid\Phi_\#r)
    +
    \int_Y H(p^y\mid r^y)\,
    d(\Phi_\#p)(y).
\]
\end{proposition}
For $\mu\in \mathcal{P}(\R^d)$ we define the entropy (i.e.\  the negative Shannon entropy)
\begin{equation}
\ent(\mu)=
\begin{cases}
     \int_{\R^d}\mu\log\mu \,dx \qquad &\text{ if }\mu\ll dx,\\
    +\infty \qquad &\text{otherwise}
\end{cases}
\end{equation}

and the Fisher information 
\begin{equation}
\label{e:fisher}
I(\mu)=
\begin{cases}
     \int_{\R^d}\abs{\nabla\log{\mu}}^2 \,d\mu, \qquad &\text{ if }\mu\ll dx\\
    +\infty \qquad &\text{otherwise}.
\end{cases}
\end{equation}
We now extend \eqref{e:entr} and define the relative entropy with respect to the reversible Wiener measure, which is an unbounded measure.
Let $\Omega=C([0,1];\R^d)$ and consider $R^\nu$ a reversible Wiener measure as defined in the last section. For $q\in\mathcal{P}(\Omega)$ with $q_0\in \mathcal{P}_2(\R^d)$ we set $\Tilde{R}^\nu(d\omega)=e^{-\pi\abs{\omega(0)}^2}R^{\nu}(d\omega)$ and define (with abuse of notation)
\begin{equation}
\label{e:entrtrick}
    H(q|R^\nu)=H(q|\Tilde{R}^\nu)-\pi\int_{\R^d}\abs{x}^2\,dq_0(x).
\end{equation}
For more insight on the definition of entropy with respect to unbounded measure, see \cite{survey,leonard_path}.
To finish the section, we recall some basic estimates regarding heat kernels and entropy, of which we provide the proof for convenience.
\begin{lemma}
\label{lem:entr}
    Let $\nu\geq  0$, $\mu\in \mathcal{P}_2(\R^d)$ and consider $\mu^{\nu}=\mu\ast \tau_\nu$. Then 
    \begin{equation}
    \label{e:secondmom}
       \int_{\R^d}\abs{x}^2\,d\mu^\nu=\int_{\R^d}\abs{x}^2\,d\mu+d\nu, 
    \end{equation}
    and
    \begin{equation}
        \label{e:wassmoll}
        W^2_2(\mu,\mu^\nu)\leq d\nu.
    \end{equation}
    Moreover,
    \begin{equation}
    \label{e:entrexp}
    -\pi \int_{\R^d}\abs{x}^2\,d\mu-\pi d\nu\leq \ent(\mu^\nu)\leq -\frac{d}{2}\log (2\pi \nu).
    \end{equation}
\end{lemma}
\begin{proof}
   We start with \eqref{e:secondmom}. Since $\mu^\nu$ is the law of a random variable $Y=X+Z$, where $X\sim \mu$ and $Z\sim \mathcal{N}(0,\nu I_d)$ are independent, we have
   \begin{equation}
       \int_{\R^d}\abs{x}^2\,d\mu^\nu=\mathbb{E}[Y^2]=\mathbb{E}[X^2]+\mathbb{E}[Z^2]=\int_{\R^d}\abs{x}^2\,d\mu+d\nu.
   \end{equation} 
\begin{comment}
Now, for $x\in\R^d$ let $\tau^x_\nu=\tau_\nu(\cdot-x)$.
   For \eqref{e:wassmoll}, consider the coupling $\gamma(dx,dy)=\tau_\nu^x(dy)\otimes \mu(dx)$. It can be checked that $\gamma\in\Pi(\mu,\mu^\nu)$, and 
   \begin{equation}
       W^2_2(\mu,\mu^\nu)=\int_{\R^d\times\R^d}\abs{x-y}^2\,d\gamma=\int_{\R^d}\int_{\R^d}\abs{y-x}^2\tau_{\nu}(y-x)\,dyd\mu(x)=M(\tau_\nu)=d\nu^2
   \end{equation}
\end{comment}
For \eqref{e:wassmoll}, consider the coupling $(X,Y)$.
We have 
\begin{equation}
    W^2_2(\mu,\mu^\nu)\leq \mathbb{E}[\abs{Y-X}^2]=\mathbb{E}[\abs{Z}^2]=d\nu.
\end{equation}
   For the lower bound of \eqref{e:entrexp}, consider the probability measure $e^{-\pi\abs{x}^2}$. Since 
    \begin{equation}
       0\leq H(\mu^\nu|e^{-\pi\abs{x}^2})=\pi \int_{\R^d}\abs{x}^2\,d\mu^\nu+\ent(\mu^\nu),
    \end{equation}
    one has 
    \begin{equation}
        \ent(\mu^\nu)\geq -\pi \int_{\R^d}\abs{x}^2\,d\mu^\nu=-\pi \int_{\R^d}\abs{x}^2\,d\mu-\pi d\nu,
    \end{equation}
    where in the last passage we use \eqref{e:secondmom}. For the upper bound, we may suppose $\nu>0$, otherwise there is nothing to prove. We notice that for all $x\in \R^d$
    \begin{equation}
        \tau_\nu(x)\leq \frac{1}{(2\pi\nu)^\frac{d}{2}},
    \end{equation} 
    thus for all $x\in \R^d$
    \begin{equation}
        \mu^\nu(x)=\int_{\R^d} \tau_\nu(x-y)d\mu(y)\leq \frac{1}{(2\pi\nu)^\frac{d}{2}}.
    \end{equation}
    This implies 
    \begin{equation}
      \log(\mu^\nu(x))\leq -\frac{d}{2}\log (2\pi \nu),
    \end{equation}
    and the conclusion follows by multiplying by $\mu^\nu(x)$ and integrating on $\R^d$. 
\end{proof}
\begin{lemma}
\label{lem:moll}
Let $\gamma\in \Pi(\rho_0,\rho_1)$. 
Consider the disintegration $\gamma(dx,dy)=\gamma_x(dy)\otimes\rho_0(dx)$. Consider for $\rho_0-$a.e.\ $x$ the measure $\gamma^\nu_x=\gamma_x\ast \tau_{\nu}$. Define $\gamma^\nu(dx,dy)=\gamma^\nu_x(dy)\otimes\rho_0(dx)$.  Then $e_1\#\gamma^\nu=\rho_1\ast\tau_\nu$,  
\begin{equation}
\label{e:wassmollcoupl}
    W_2^2(\gamma^\nu,\gamma)\leq d\nu 
\end{equation}
and 
\begin{equation}
\label{e:entrgamma}
    \ent (\gamma^\nu)\leq \ent(\rho_0)-\frac{d}{2}\log(2\pi\nu).
\end{equation}
\end{lemma}
\begin{proof}
\begin{comment}
    For any $\varphi\in C_0(\R^d)$, by using the property $\int_{\R^d} f(x)d(\mu\ast h)(x)=\int_{\R^d} f\ast h(x)d\mu(x)$,
    \begin{multline}
        \int_{\R^d\times \R^d}\varphi(y)d\gamma^\nu(x,y)=
        \int_{\R^d\times \R^d}\varphi(y)d(\gamma_x\ast\tau_\nu)(y)d\rho_0(x)=\int_{\R^d\times \R^d}\varphi\ast\tau_\nu(y)d\gamma_x(y)d\rho_0(x)\\=\int_{\R^d\times \R^d}\varphi\ast\tau_\nu(y)d\gamma(x,y)=\int_{\R^d}\varphi\ast\tau_\nu(y)d\rho_1(y)=\int_{\R^d}\varphi(y)d(\rho_1\ast\tau_\nu)(y),
    \end{multline}
    which means that $e_1\#\gamma^\nu=\rho_1\ast\tau_\nu$.
\end{comment}
 Consider random variables $(X,Y)\sim \gamma$ and $Z\sim \tau_\nu$ independent of $(X,Y)$. Consider the vector $(X,Y+Z)$. By conditioning one finds that $(X,Y+Z)\sim \gamma^\nu$. In particular this implies $e_1\#\gamma^\nu=\operatorname{Law}(Y+Z)=\rho_1\ast\tau_\nu$. Let us now prove \eqref{e:wassmollcoupl}. Choosing the coupling $((X,Y),(X,Y+Z))$, we find 
    \begin{equation}
        W_2^2(\gamma^\nu,\gamma)\leq \mathbb{E}\left[\abs{(X,Y)-(X,Y+Z)}^2\right]=\mathbb{E}\left[\abs{Z}^2\right]\leq d\nu.      
    \end{equation}
    Finally, for \eqref{e:entrgamma}, by Proposition \ref{prop:entrdisint} applied to $e_0$
     \begin{equation}
    \ent(\gamma^\nu)=\ent(\rho_0)+\int_{\R^d}\ent(\gamma^\nu_x)d\rho_0(x)\overset{\eqref{e:entrexp}}{\leq}\ent(\rho_0)-\frac{d}{2}\log(2\pi\nu).  
     \end{equation}
    \end{proof}
    \subsection{Hypotheses and Functionals}
    \label{sec:hyp}
Let us now specify the precise formulation of the problem and the set of hypotheses we make on the data.
\begin{hypothesis}[Initial condition]
\label{hyp:initial}
We fix an initial probability measure
$\rho_0\in\mathcal P_2(\mathbb R^d)$ satisfying
\[
\operatorname{Ent}(\rho_0)<+\infty.
\]
\end{hypothesis}

We consider two alternative sets of hypotheses on the
interaction functionals
\[
F,G:\mathcal P_2(\mathbb R^d)
\longrightarrow [0,+\infty].
\]

\begin{hypothesis}[Mean-Field Games]
\label{hyp:mfg}
Similarly to \cite{benamou2019entropy}, we suppose that 
\begin{equation}
    F=F_1+F_2, \qquad G=G_1+G_2,
\end{equation}
where for some convex continuous $L_F,L_G:[0,+\infty)\rightarrow\R^+$ with $L_F(0)=L_G(0)=0$
\begin{equation}
    F_1(\rho)=\begin{cases}
        \int_{\R^d}L_F(\rho(x))\,dx \quad & \text{ if } \rho\ll dx\\
        +\infty \quad & \text{ otherwise, }
    \end{cases}
\end{equation}
\begin{equation}
    G_1(\rho)=\begin{cases}
        \int_{\R^d}L_G(\rho(x))\,dx \quad & \text{ if } \rho\ll dx\\
        +\infty \quad & \text{ otherwise, }
    \end{cases}
\end{equation}
and $F_2,G_2:\mathcal{P}_2(\R^d)\rightarrow [0,+\infty)$ are continuous with respect to $W_2$ and narrowly lower semicontinuous. We furthermore assume, coherently with \cite{cardaliaguetweak,cardaliaguetfo,variationalMFG}\footnote{ In \cite{cardaliaguetweak,cardaliaguetfo,variationalMFG}, the coupling $f$ in \eqref{e:thesystem} is assumed to be completely local, and they ask for conditions of the form \begin{equation}
    \frac{1}{C}\abs{\rho}^{q-1}-C\leq f(\rho)\leq C\abs{\rho}^{q-1}+C
\end{equation}
for some $C>0$ and $q>1$. Since in our setting $L_F$ here is an antiderivative of the local part of $f$, this condition translates to \eqref{e:lfass}. A similar reasoning holds for \eqref{e:lgass}.}, that there exist $q_F,q_G>1$ and $C_F,C_G>0$ such that for all $\rho>0$
\begin{equation}
\label{e:lfass}
    \frac{1}{C_F}\abs{\rho}^{q_F}-C_F\rho\leq L_F(\rho)\leq C_F(\abs{\rho}^{q_F}+\rho),
\end{equation}
\begin{equation}
\label{e:lgass}
    \frac{1}{C_G}\abs{\rho}^{q_G}-C_G\rho\leq L_G(\rho)\leq C_G(\abs{\rho}^{q_G}+\rho).
\end{equation}
\end{hypothesis}

\begin{hypothesis}[Entropic optimal transport]
\label{hyp:eot}
Alternatively, we suppose that,
given
$\rho_1\in\mathcal P_2(\mathbb R^d)$
with $\operatorname{Ent}(\rho_1)<+\infty$,
\[
F\equiv0, \qquad G(\rho)=\iota_{\rho_1}(\rho)
=
\begin{cases}
0, & \text{if }\rho=\rho_1,\\
+\infty, & \text{otherwise}.
\end{cases}
\]
This, paired with \ref{hyp:initial}, is the standard setting of entropic optimal transport.
\end{hypothesis}

Notice that \eqref{e:lfass} and \eqref{e:lgass} immediately imply that for any $\rho\in \mathcal{P}(\R^d)$ with $\rho\ll dx$
\begin{equation}
\label{e:lpf}
    \frac{1}{C_F}\norm{\rho}_{L^{q_F}(\R^d)}^{q_F}-C_F\leq F_1(\rho)\leq C_F( \norm{\rho}_{L^{q_F}(\R^d)}^{q_F}+1),
\end{equation}
\begin{equation}
\label{e:lpg}
    \frac{1}{C_G}\norm{\rho}_{L^{q_G}(\R^d)}^{q_G}-C_G\leq G_1(\rho)\leq C_G( \norm{\rho}_{L^{q_G}(\R^d)}^{q_G}+1).
\end{equation}
Moreover, a standard argument employing the continuity of $L_F$, paired with the right hand side inequality in \eqref{e:lpf} and the Dominated Convergence Theorem gives that $F_1(\rho)$ is continuous on $L^{q_F}(\R^d)$. Similarly, $G_1(\rho)$ is continuous on $L^{q_G}(\R^d)$.
For $\rho\in C([0,1],\mathcal{P}_2(\R^d))$ we set
\begin{equation}
    C(\rho)=\int_0^1F(\rho_t)dt+G(\rho_1).
\end{equation}
We also set for convenience 
\begin{equation}
    C_1(\rho)=\int_0^1F_1(\rho_t)dt+G_1(\rho_1),
\end{equation}
\begin{equation}
    C_2(\rho)=\int_0^1F_2(\rho_t)dt+G_2(\rho_1).
\end{equation}

For any $\rho\in C^0([0,1],\mathcal{P}_2(\R^d))$ and $v\in L^{2}([0,1]\times \R^d,dt \otimes d\rho_t)$ we introduce the (possibly infinite) kinetic energy
\begin{equation}
    \mathcal{A}(\rho,v)=\frac{1}{2}\int_{0}^1\int_{\R^d}\abs{v}^2d\rho_t dt 
\end{equation} 
and (with abuse of notation)
\begin{equation}
    \mathcal{A}(\rho)=\inf_{v}\left\{\mathcal{A}(\rho,v): \partial_t\rho+\dive(\rho v)=0 \right\},
\end{equation}
where the continuity equation has to be understood distributionally (see \cite[Chapter 8]{AGS}).
We say that $\rho\in AC^2([0,1],\mathcal{P}_2(\R^d))$ if $\mathcal{A}(\rho)<+\infty$ (see \cite[Theorem 8.3.1]{AGS}). For $\nu\geq 0$ and $\rho \in C^{0}([0,1],\mathcal{P}_2(\R^d))$ we define
\begin{equation}
\label{e:Enu}
    \mathcal{E}_{\nu}(\rho)=\mathcal{A}(\rho)+\frac{\nu^2}{8}\int_0^1I(\rho_t)dt+ \frac{\nu}{2}(\ent(\rho_1)-\ent(\rho_0)).
\end{equation}
Moreover, for $\rho_0 \in \mathcal{P}_2({\R^d})$ let 
\begin{equation}
    \iota_{\rho_0}(\rho)=\begin{cases}
        0\quad &\text{if } \rho|_{t=0}=\rho_0,\\
        +\infty &\text{otherwise }.
    \end{cases}
\end{equation}
For $\nu\geq 0$ and $\rho_0\in \mathcal{P}_2({\R^d})$ we define the Eulerian energy as $\mathcal{F}_\nu:C^{0}([0,1],\mathcal{P}_2(\R^d))\rightarrow\R\cup \{+\infty\}$ as
\begin{equation}
\label{e:eulene}
    \mathcal{F}_\nu(\rho)=\mathcal{E}_{\nu}(\rho)+\iota_{\rho_0}(\rho)+C(\rho),
\end{equation}
so that the Eulerian minimization problem reads 
\begin{equation}
\label{e:probeul}
    \inf_{\rho\in C^0([0,1],\mathcal{P}_2(\R^d))}\mathcal{F}_\nu (\rho).
\end{equation}
We now introduce the Lagrangian functional. For $Q\in\mathcal{P}(\Omega)$ we let \begin{equation}
    \mathcal{K}(Q)=\int_\Omega k(\omega)dQ(\omega)
\end{equation}
where
    \begin{equation}
        k(\omega)=\begin{cases}
             \int_0^1\frac{1}{2}\abs{\dot{\omega}(t)}^2dt& \text{ if } \omega\in AC^2([0,1];\R^d),\\
    +\infty & \text{ otherwise }.
        \end{cases}
    \end{equation}
    $\mathcal{K}$ is narrowly lower semicontinuous, see for instance \cite{brenieraction}.
We moreover define (with abuse of notation) $C(Q)=C((Q_t)_{t\in [0,1]})$ and $\iota_{\rho_0}(Q)=\iota_{\rho_0}((Q_t)_{t\in [0,1]})$.
For $\nu >0$ and $\rho_0\in \mathcal{P}_2({\R^d})$ we now define the Lagrangian functional $\mathcal{L}_\nu:\mathcal{P}(\Omega)\rightarrow\R\cup \{+\infty\}$ as
\begin{equation}
    \mathcal{L}_\nu(Q)=\nu H(Q|R^{\nu})+C(Q)+\iota_{\rho_0}(Q),
\end{equation}
where $R^\nu$ is the reversible Wiener measure, so that the Lagrangian minimization problem reads
\begin{equation}
\label{e:problagr}
    \inf_{Q\in \mathcal{P}(\Omega)}\mathcal{L}_\nu(Q).
\end{equation}
Define moreover the Lagrangian limit problem $\mathcal{L}_0:\mathcal{P}(\Omega)\rightarrow\R\cup \{+\infty\}$ as
\begin{equation}
    \mathcal{L}_0(Q)=\mathcal{K}(Q)+C(Q)+\iota_{\rho_0}(Q).
\end{equation} 
\section{Eulerian Limit}
\label{sec:eul}
Let us start by assessing existence.
\begin{proposition}
\label{prop:existenceeul}
 Let $\nu\geq 0$. Suppose \ref{hyp:initial} and either \ref{hyp:mfg} or \ref{hyp:eot} hold. Suppose moreover there exists $\overline{\rho}\in C^0([0,1];\mathcal{P}_2(\R^d))$ such that $\mathcal{F}_\nu(\overline{\rho})<+\infty$. 
     Then for all $\nu\geq 0$ the problem 
\eqref{e:probeul}
    admits a minimizer.
    \end{proposition}
    \begin{proof}
        In the case of \ref{hyp:eot}, this is just the existence of an optimal (entropic) transport plan \cite{survey}. In the case of \ref{hyp:mfg}, see for instance \cite{variationalMFG,lavenantmfg,cardaliaguetweak,cardaliaguetfo}.
    \end{proof}
    As usual, we want to pair the $\Gamma$-convergence result with an equicompactness statement, to deduce convergence of minimizers.
   In our case this role is played by the following quite standard result for the action functional.
\begin{proposition}
\label{prop:compeul}
Let $\rho_0\in\mathcal{P}_2(\R^d)$. Let $\{\rho^n\}_{n\in\N}\subset C([0,1],\mathcal{P}_2({\R^d}))$ such that $ \mathcal{A}(\rho^n)+\iota_{\rho_0}(\rho^n)\leq C$. Then there exists a subsequence $\rho^{n_k}\rightarrow \rho\in C([0,1],\mathcal{P}_2({\R^d}))$ narrowly.
\end{proposition}
\begin{proof}
    Let $t,s\in[0,1]$. Using \eqref{e:bb} (after rescaling the time interval between $s$ and $t$ to $[0,1]$), 
    \begin{equation}
    \label{e:equiconeu}
        W_2^2(\rho_t^n,\rho_s^n)\leq  2\mathcal{A}(\rho^n)\abs{t-s}\lesssim \abs{t-s}. 
    \end{equation}
     Moreover, 
  \begin{equation}
      \int_{\R^d}\abs{x}^2d\rho_t^n=W_2^2(\rho_t^n,\delta_0)\lesssim W_2^2(\rho_0,\delta_0)+W_2^2(\rho_0,\rho_t^n)\overset{\eqref{e:equiconeu}}{\lesssim}1,  
  \end{equation}
  hence Ascoli-Arzelà theorem on $C^0([0,1],\mathcal{P}(\R^d))$ with the metric \eqref{e:narrow} (see \cite[Proposition 3.3.1]{AGS}) gives the existence of a converging subsequence. Finally we are left to prove that $\rho\in C^0([0,1],\mathcal{P}_2(\R^d))$, but this is a direct consequence of \eqref{e:equiconeu} since $W_2^2(\cdot,\cdot)$ is lower semicontinuous with respect to narrow convergence \cite[Proposition 7.1.3]{AGS}.
\end{proof}
\begin{remark}
\label{rem:comp}
    It is in general false that sequences with equibounded action functional admit a subsequence uniformly converging in $C([0,1],\mathcal{P}_2({\R^d}))$. Indeed,  consider for instance for $d=1$ 
the curve
\[
    \rho_t^n
    = \left(1-\frac{1}{n^2}\right)\delta_0
      + \frac{1}{n^2}\delta_{nt},
    \qquad t\in[0,1].
\]
The curve is admissible for the continuity equation:
the mass $1/n^2$ moves at constant velocity $n$,
while the remaining mass stays at the origin. Hence,
\[
    \mathcal A(\rho^n)
    \le \frac12\int_0^1
        \frac{1}{n^2}n^2\,dt
    =\frac12.
\]
For every $t\in[0,1]$,
    $\rho_t^n \rightarrow \delta_0$ narrowly, however
\[
    W_2^2(\rho_t^n,\delta_0)
    =\int_{\mathbb R}|x|^2\,d\rho_t^n(x)
    =\frac{1}{n^2}(nt)^2=t^2.
\]
A similar counterexample can be constructed also if we assume equiboundedness of the full energy \eqref{e:eulene}.
\end{remark}
\subsection{Eulerian $\Gamma$-Limit}
In this section we prove the following theorem. 
\begin{theorem}
\label{th:maineul}
Suppose \ref{hyp:initial} and either \ref{hyp:mfg} or \ref{hyp:eot} hold. Then, $\mathcal{F}_{\nu}\overset{\Gamma}{\rightarrow}\mathcal{F}_0$ with respect to narrow convergence \eqref{e:narrow} in $C([0,1],\mathcal{P}({\R^d}))$.
\end{theorem}
The main ingredient is the following lemma. 
\begin{lemma}
\label{lem:limsup}
Let $\rho\in C([0,1],\mathcal{P}_2({\R^d}))$ with $\mathcal{E}_0(\rho)<\infty$. Let $\tau_\nu$ be the heat kernel. 
Let $\rho^\nu_t=\rho_t\ast \tau_{\nu t}$.
    Then
    \begin{equation}
        \mathcal{E}_\nu(\rho^\nu)\leq \mathcal{A}(\rho).
    \end{equation}
\end{lemma}
\begin{proof}
    The proof is an easier case of \cite[Lemma 4.1]{baradatm}. We provide the main steps. Consider the velocity field $v$ of $\rho$. We define 
    \begin{equation}
        \hat{v}^{\nu}_t=\frac{(\rho_tv_t)\ast\tau_{\nu t}}{\rho_t^{\nu}}.
    \end{equation}
 By a direct computation, $(\rho^{\nu},\hat{v}^{\nu})$ solves 
 \begin{equation}
 \partial_t\rho_t^{\nu}+\dive(\rho_t^{\nu}\hat{v}^{\nu}_t)=\frac{\nu}{2}\Delta\rho_t^{\nu}.
 \end{equation}
 By defining ${v}_t^{\nu}=\hat{v}_t^{\nu}-\frac{\nu}{2}\nabla(\log\rho_t^{\nu})$, the couple $(\rho^{\nu},v^{\nu})$ solves 
 \begin{equation}
     \partial_t\rho_t^{\nu}+\dive(\rho_t^{\nu}v^{\nu}_t)=0.
 \end{equation}
 By definition 
\begin{equation}
\label{e:comp}
\mathcal{A}(\rho^\nu)\leq\mathcal{A}(\rho^\nu,v^{\nu}),     
\end{equation}
while by Jensen's inequality 
\begin{equation}
\label{e:jens}
   \frac{1}{2} \int_0^1\int_{\R^d}\abs{v_t^\nu+\frac{\nu}{2}\nabla(\log\rho_t^\nu)}^2\,d\rho_t^\nu dt=\mathcal{A}(\rho^\nu,\hat{v}^{\nu})\leq\mathcal{A}(\rho,v)=\mathcal{A}(\rho). 
\end{equation}
Combining the two inequalities \eqref{e:comp} and \eqref{e:jens}, and expanding the square in \eqref{e:jens} (with the same computation as in \cite[Lemma 4.1]{baradatm}), we find the result.

\end{proof}

\begin{comment}
We also need that the forcing terms $F,G$ are compatible with the topology of the $\Gamma-$convergence.
\begin{lemma}
\label{lem:cont}The following continuity properties hold:
\begin{enumerate}
    \item \label{it:1} The functionals $\rho\rightarrow\int_{\R^d}L_F(\rho(x))\,dx$ and $\rho\rightarrow\int_{\R^d}L_G(\rho(x))\,dx$ are weakly semicontinuous on $L^1(\R^d)$. \AC{this to be combined with boundedness of the energy: Delavallé-Pussin + Dunford-Pettis?}
        \item \label{it:2} Let $\rho_n\rightarrow \rho$ in $C([0,T];L^1({\R^d}))$ with $\sup_n C_1(\rho_n)<+\infty$. Then $C_1(\rho_n)\rightarrow C(\rho)$. \AC{Maybe need more assumptions on $L_F,L_G$, I think we need $L^q$ topology and $f$ bounded from above and below by $m^q$}
    \item \label{it:3} The functional $C_2$ is continuous on $C([0,T];\mathcal{P}_2({\R^d}))$. 
\end{enumerate}
\end{lemma}
\begin{proof}
   For \ref{it:3}, it suffices to notice that $F_2$ and $G_2$ are continuous with respect to $W_2$.
\end{proof}
\end{comment}
\begin{proof}[Proof of Theorem \ref{th:maineul}]
    The proof under \ref{hyp:eot} is a direct extension of \cite[Theorem 2.3]{baradatm} to $\R^d$, thus we skip it. Let us start by $\Gamma$-$\liminf$ inequality. Let $\nu_n\rightarrow 0$  and $\rho^{n}\rightarrow \rho$ narrowly in $C([0,1],\mathcal{P}({\R^d}))$. Up to choosing a subsequence which realizes the $\liminf$, we may suppose that $\mathcal{F}_{\nu_n}(\rho^n)\leq C$, otherwise there is nothing to prove. We first show that $\mathcal A(\rho^n)$ is bounded. By \eqref{e:bb}, $\int_{\R^d}\abs{x}^2d\rho^n_1\leq 2\int_{\R^d}\abs{x}^2d\rho_0+4\mathcal A(\rho^n)$, and $\ent(\mu)\geq-\pi\int\abs{x}^2d\mu$ for every $\mu\in\mathcal P_2(\R^d)$ (see Lemma \ref{lem:entr}). Since $F,G\geq 0$ and $I\geq 0$,
\begin{equation}
C\geq\mathcal F_{\nu_n}(\rho^n)\geq(1-2\pi\nu_n)\mathcal A(\rho^n)-\frac{\nu_n}{2}\ent(\rho_0)-\pi\nu_n\int_{\R^d}\abs{x}^2d\rho_0,
\end{equation}
so that $\mathcal A(\rho^n)\leq C$ for $n$ large. Arguing as in Proposition \ref{prop:compeul}, it then holds 
    \begin{equation}
        \int_{\R^d} \abs{x}^2d\rho_1^{n}\leq C.
    \end{equation}
    Hence, again by the lower bound $\ent\geq-\pi\int\abs{x}^2$,
    \begin{equation}
        \nu_n \ent(\rho^{n}_1)\geq  -\nu_n \pi \int_{\R^d} \abs{x}^2d\rho_1^{n}\rightarrow 0.
    \end{equation} 
    By the positivity of $I$ (see \eqref{e:fisher}) and the narrow lower semicontinuity of $\mathcal A$ (see e.g.\ \cite[Section 5.3]{santambrogio2015optimal}), looking at \eqref{e:Enu} we deduce \begin{equation}
    \liminf_{n}\mathcal{E}_{\nu_n}(\rho^n)\geq \mathcal{E}_{0}(\rho). 
    \end{equation}
    To conclude the $\Gamma$-$\liminf$ inequality, we prove that $C(\rho)$ is lower-semicontinuous with respect to the narrow convergence on $C([0,1],\mathcal{P}({\R^d}))$. For the term $C_2(\rho)$, this is immediate since both $F_2$ and $G_2$ are narrowly lower semicontinuous on $\mathcal{P}(\R^d)$ and convergence is uniform in time. 
    Notice now that $L_F$ is convex, and that the left hand side inequality in \eqref{e:lfass} implies that $L_F$ is superlinear, that is \begin{equation}
        \lim_{\rho\rightarrow +\infty}\frac{L_F(\rho)}{\rho}=+\infty.
    \end{equation}
    Hence, by \cite[Proposition 7.7]{santambrogio2015optimal}, $F_1$ is narrowly lower semicontinuous on $\mathcal{P}(\R^d)$. Similarly, \eqref{e:lgass} implies that $G_1$ is narrowly lower semicontinuous on $\mathcal{P}(\R^d)$. Using Fatou's Lemma 
    \begin{equation}
        \liminf_n\int_0^1F_1(\rho_t^n)\,dt\geq \int_0^1\liminf_n F_1(\rho_t^n)\,dt\geq \int_0^1F_1(\rho_t)\,dt,
    \end{equation}
    which concludes the proof of the lower semicontinuity of $C_1$, and thus the $\Gamma$-$\liminf$ inequality.

   For the $\Gamma$-$\limsup$ inequality, consider $\rho\in AC^2([0,1],\mathcal{P}_2({\R^d}))$ such that $\mathcal F_0(\rho)<+\infty$. Consider for $\nu_n\rightarrow 0$ the heat kernel $\tau_{\nu_n}$ and define $\rho^n_t=\rho_t\ast\tau_{\nu_n t}$. A direct computation using \eqref{e:wassmoll} gives $\rho^n\rightarrow \rho $ in $C^0([0,1],\mathcal{P}_2(\R^d))$.  
    By Lemma \ref{lem:limsup}, we have 
    \begin{equation}
         \limsup_{n}\mathcal{E}_{\nu_n}(\rho^n)\leq \mathcal{E}_{0}(\rho).
    \end{equation}
    To conclude, we just need to prove that $C(\rho^n)\rightarrow C(\rho)$. The fact that $C_2(\rho^n)\rightarrow C_2(\rho)$ is an immediate consequence of the properties of $F_2,G_2$. For $C_1$, notice that $\mathcal{F}_0(\rho)<+\infty$ combined with the left hand side inequality in \eqref{e:lpf} gives that $\norm{\rho_t}_{L^{q_F}}^{q_F}\in L^{1}(0,1)$, so that $\rho_t\in L^{q_F}(\R^d)$ for a.e.\ $t\in [0,1]$. Similarly, \eqref{e:lpg} gives $\rho_{1}\in L^{q_G}(\R^d)$. Since $\rho^n$ is defined as convolution, we thus have $\rho^n_t\rightarrow\rho_t$ in $ L^{q_F}(\R^d)$ for a.e.\ $t\in [0,1]$ and $\rho^n_1\rightarrow\rho_1$ in $ L^{q_G}(\R^d)$. Hence, by the discussion after \eqref{e:lpg}, $G_{1}(\rho^n_1)\rightarrow G_{1}(\rho_1)$ and, for a.e.\ $t\in [0,1]$, $F_1(\rho_t^n)\rightarrow F_1(\rho_t)$. Since by \eqref{e:lpf} and properties of the convolution
    \begin{equation}
        F_1(\rho^n_t)\leq C_F(\norm{\rho_t^n}_{L^{q_F}(\R^d)}^{q_F}+1)\leq C_F(\norm{\rho_t}_{L^{q_F}(\R^d)}^{q_F}+1)\in L^{1}(0,1), 
    \end{equation}
    an application of the Dominated Convergence Theorem concludes the proof of the continuity of $C_1$, and thus finishes the proof of the $\Gamma$-$\limsup$ inequality.
\end{proof}
\begin{remark}
 Looking at the proof of Theorem \ref{th:maineul}, one realizes that the same result holds also with respect to the topology given by \eqref{e:uniform}. We choose to state the result this way since it matches the compactness result given by Proposition \ref{prop:compeul}.
\end{remark}
\subsection{First Order Energy Expansion}
\label{sec:expansion}

The goal of this section is to derive, in the case of \ref{hyp:mfg}, an energy expansion in the spirit of the works \cite{giglitam,confortitam,legerexpansion} for \ref{hyp:eot}.
While not explicitly referring to notions of higher order convergence, we follow the notation of \cite{higherordergammaconv}. In this section we use the shorthand notation $\mathcal{F}=\mathcal{F}_0$. 
 Moreover, set 
\begin{equation}
    \mathcal{U}_0=\argmin_{\rho\in AC^2([0,1],\mathcal{P}_2(\R^d))} \mathcal{F},
\end{equation}
that is the set of minimizers of the first order MFG. Since the non-convex terms $F_2$ and $G_2$ are present, no uniqueness is expected for the problem. 
Define \begin{equation}
\mathcal{F}^1(\rho)=
        \frac{1}{2}(\ent(\rho_1)-\ent(\rho_0)),
\qquad
\mathcal{F}^2(\rho)=\frac{1}{8}\int_0^1I(\rho_t)\,dt
\end{equation}
and
\begin{equation}
    m_1=\inf_{\rho\in \mathcal{U}_0}\mathcal{F}^1(\rho), \qquad \mathcal{U}_1=\argmin_{\rho\in \mathcal{U}_0}\mathcal{F}^1(\rho).
\end{equation}
\begin{theorem}
\label{th:energyexp}
Suppose \ref{hyp:initial} and \ref{hyp:mfg} hold. If  $m_1<+\infty$ and $\mathcal{U}_1\cap \{\mathcal{F}^2(\rho)<+\infty\}\neq \emptyset$, then 
\begin{equation}
\label{e:energyexp}
    m_{\nu}=m_0+\nu m_1+o(\nu).
\end{equation}
Moreover, if $\nu_n\rightarrow 0$ and $\{\rho_{n}\}_{n\in\N}$ is a sequence of minimizers $\rho_n\in \argmin \mathcal{F}_{\nu_n}$, then any cluster point $\rho$ of $\{\rho_{n}\}_{n\in\N}$ belongs to $\mathcal{U}_1$. 
\end{theorem}
\begin{proof}
Let us address first \eqref{e:energyexp}.
Take $\overline{\rho}\in \mathcal{U}_1\cap \{\mathcal{F}^2(\rho)<+\infty\}$. Since $m_\nu\leq \mathcal{F}_\nu(\overline{\rho})$ it holds 
\begin{equation}
\label{e:upperb}
    m_\nu\leq \mathcal{F}(\overline{\rho})+\nu \mathcal{F}^1(\overline{\rho})+\nu^2\mathcal{F}^2(\overline{\rho})=m_0+\nu m_1+\nu^2\mathcal{F}^2(\overline{\rho}),
\end{equation}
so that 
\begin{equation}
\limsup_\nu \frac{m_\nu-m_0}{\nu}\leq m_1.     
\end{equation}
We are left to prove 
\begin{equation}
   \liminf_\nu \frac{m_\nu-m_0}{\nu}\geq m_1.  
\end{equation}
Take $\nu_n\rightarrow 0$ which realizes the liminf. Consider $\{\rho_{n}\}_{n\in\N}$ a sequence of minimizers $\rho_n\in \argmin \mathcal{F}_{\nu_n}$. By the bound on $\mathcal A(\rho_n)$ obtained in the proof of Theorem \ref{th:maineul} and Proposition \ref{prop:compeul}, up to a subsequence, there exists $\rho$ such that $\rho_n\rightarrow \rho$ narrowly, and by the fundamental Theorem of $\Gamma$-convergence $\rho\in \mathcal{U}_0$. 
This means that, since $\mathcal{F}(\rho_n)\geq m_0$ and $\mathcal{F}^2(\rho_n)\geq 0$, 
\begin{multline}
    \liminf_\nu \frac{m_{\nu}-m_0}{\nu}=\lim_n \frac{m_{\nu_n}-m_0}{\nu_n}=\lim_n\left[ \frac{\mathcal{F}(\rho_n)-m_0}{\nu_n}+\mathcal{F}^1(\rho_n)+\nu_n\mathcal{F}^2(\rho_n)\right]\\ \geq \liminf_n\mathcal{F}^1(\rho_n).
\end{multline}
Finally, by lower semi-continuity of $\mathcal{F}^1$ ($\ent(\rho_0)$ is a constant and $\ent$ is lower semicontinuous along sequences whose second moments are bounded) and since $\rho\in \mathcal{U}_0$, 
\begin{equation}
    \liminf_n\mathcal{F}^1(\rho_n) \geq \mathcal{F}^1(\rho)\geq m_1.
\end{equation}
This concludes the proof of \eqref{e:energyexp}.
Repeating the above computations for any $\nu_n\rightarrow 0$ and $\{\rho_n\}_{n\in\N}$ sequence of minimizers for $\mathcal{F}_{\nu_n}$, we get that, for a cluster point $\rho$ of $\{\rho_n\}_{n\in\N}$, necessarily $\mathcal{F}^1(\rho)=m_1$.
\end{proof}
\begin{remark}
\label{rem:exp}
 Unlike \cite{giglitam,confortitam,legerexpansion}, the energy expansion \eqref{e:energyexp} is only of order $1$. Our argument, which is inspired by  \cite{giglitam,confortitam,legerexpansion}, does not translate directly to a second order expansion, since in the MFG problem, differently from the Schrödinger problem,  the final marginal is in general not fixed. Indeed, we see that, while a second order limsup inequality would hold with an analogous proof as above, for the liminf inequality one finds 
    \begin{equation}
        \liminf_\nu \frac{m_{\nu}-\nu m_1-m_0}{\nu^2}= \lim_n \left[\frac{\mathcal{F}(\rho_n)-m_0}{\nu_n^2}+\frac{\mathcal{F}^1(\rho_n)-m_1}{\nu_n}+\mathcal{F}^2(\rho_n)\right].
    \end{equation}
    Looking at the right hand side, if the final marginal is not fixed, the second term could be negative preventing to deduce the inequality.
\end{remark}

\begin{figure}[h]
    \centering\includegraphics{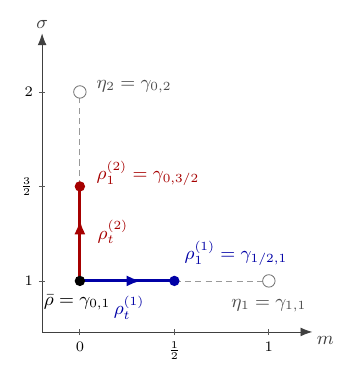}
    \caption{The two minimizers, plotted with respect to their mean $m$ and standard deviation $\sigma$.}
    \label{fig:twominima}
\end{figure}

We finish the section by showing an example of selection of minimizers.
Let \(d=1\), and call $
\gamma_{m,\sigma}=\mathcal N(m,\sigma^2)$
the Gaussian probability measure with mean \(m\) and standard deviation \(\sigma\). Take the initial condition
$\rho_0=\gamma_{0,1}$, $F=0$ and
\[
G(\mu)=G_2(\mu)
=
\frac12\min\left\{
W_2^2(\mu,\eta_1),\,W_2^2(\mu,\eta_2)
\right\},
\]
where
$
\eta_1=\gamma_{1,1},
\eta_2=\gamma_{0,2}.
$
These assumptions technically do not fit the set of hypotheses \ref{hyp:mfg} (we cannot choose $F_1,G_1=0$). However, Theorem \ref{th:energyexp} applies also in this setting: its proof only uses the existence of minimizers of $\mathcal F_\nu$, Proposition \ref{prop:compeul} and Theorem \ref{th:maineul}, which hold here since $F=0$ and $G=G_2$ is $W_2$-continuous and narrowly lower semicontinuous.
We recall that for one dimensional gaussian measures
\[
W_2^2(\gamma_{m,\sigma},\gamma_{\tilde m,\tilde\sigma})
=
|m-\tilde m|^2+|\sigma-\tilde\sigma|^2,
\]
thus
\[
W_2^2(\rho_0,\eta_1)
=
W_2^2(\rho_0,\eta_2)
=
1.
\]
We have that 
\[
m_0
=\inf_{\rho_{|t=0}=\rho_0}
\left\{
\mathcal A(\rho)+G(\rho_1)
\right\}=
\min_{j=1,2}\inf_{\mu\in\mathcal P_2(\mathbb R)}
\frac12\left[
W_2^2(\rho_0,\mu)+W_2^2(\mu,\eta_j)
\right].
\]
For $j=1,2$, by the triangle inequality,
\[
W_2(\rho_0,\mu)+W_2(\mu,\eta_j)\ge W_2(\rho_0,\eta_j)=1,
\]
 hence by Young's inequality
 \begin{equation}
\label{e:young}
\frac12(W_2(\rho_0,\mu)^2+W_2(\mu,\eta_j)^2)
\ge
\frac14(W_2(\rho_0,\mu)+W_2(\mu,\eta_j))^2
\ge
\frac14.
\end{equation}
Equality holds in \eqref{e:young} if and only if
\begin{equation}
    W_2(\rho_0,\mu)=W_2(\mu,\eta_j)\quad\text{and}\quad W_2(\rho_0,\mu)+W_2(\mu,\eta_j)=1,
\end{equation} that is \(\mu\) is the midpoint of the unique Wasserstein geodesic between \(\rho_0\) and \(\eta_j\). Thus
$
m_0=\frac14$
and $\mathcal{F}$ admits two minimizers, which are explicitly computable since $\rho_0,\eta_1,\eta_2$ are Gaussian:
\[
\rho_t^{(1)}=\gamma_{t/2,1},
\qquad
\rho_t^{(2)}=\gamma_{0,\,1+t/2},
\]

see Figure \ref{fig:twominima}.
Recall that for a one-dimensional Gaussian,
\[
\operatorname{Ent}(\gamma_{m,\sigma})
=
-\frac12\log(2\pi e\sigma^2),
\qquad
I(\gamma_{m,\sigma})=\frac1{\sigma^2}.
\]
Computing
\[
\int_0^1 I(\rho_t^{(1)})\,dt=1, \qquad \int_0^1 I(\rho_t^{(2)})\,dt=
\frac23,
\]
and 
\[
\mathcal F^1(\rho^{(1)})=0, \qquad \mathcal F^1(\rho^{(2)})
=
-\frac12\log\frac32<0.
\]
It follows that
\[
m_1=-\frac12\log\frac32,
\qquad
\mathcal U_1=\{\rho^{(2)}\}.
\]
By Theorem \ref{th:energyexp} we can conclude that 
\[
m_\nu
=
\frac14
-\frac{\nu}{2}\log\frac32
+o(\nu),
\]
and every cluster point of minimizers of $(\mathcal{F}_\nu)_\nu$ is \(\rho^{(2)}\).

\section{Lagrangian limit}
\label{sec:lagr}
In this section we prove the following theorem.
\begin{theorem}
\label{th:mainlagr}
    Suppose \ref{hyp:initial} and either \ref{hyp:mfg} or \ref{hyp:eot} hold. Then $\mathcal{L}_\nu\overset{\Gamma}{\rightarrow} \mathcal{L}_0$ with respect to narrow convergence on $\mathcal{P}(\Omega)$. \end{theorem}
    \begin{proof}
      The Theorem follows by combining Proposition \ref{prop:lower} and Proposition \ref{prop:upper}.  
    \end{proof}
We recall a definition and a lemma from \cite{baradatm}.
  \begin{definition}
  \label{def:kn}
      Let $\omega\in C([0,1],\R^d)$ and $N\in \N$. Set $h=1/N$ and define for $n=0\dots N$ the times $t_n=nh$. We consider the approximation of the kinetic action
      \begin{equation}
          k_N(\omega)=\min(\Tilde{k}_N(\omega),N)
      \end{equation}
      where 
      \begin{equation}
          \Tilde{k}_N(\omega)=\frac{1}{2h}\sum_{n=0}^{N-1}\abs{\omega(t_{n+1})-\omega(t_{n})}^2.
      \end{equation}
  \end{definition}
  The following lemma is proved in \cite[Lemma 3.3]{baradatm}.
  \begin{lemma}
  \label{lem:baradatkin}
      For $N\in \N$, the functional $k_N$ is continuous and bounded on $C([0,1],\R^d)$. For every $\omega\in C([0,1],\R^d)$ 
      \begin{equation}
          k_N(\omega)\rightarrow
          \begin{cases}
              k(\omega) &\quad \text{ if } \omega\in AC^2([0,1];\R^d),\\
              +\infty &\quad \text{ otherwise. }
          \end{cases}
      \end{equation}
      Moreover, for any $\alpha\in (0,1)$ and $\nu\in (0,1)$ there holds 
      \begin{equation}
      \label{e:knest}
          \int \exp\left(\frac{\alpha}{\nu}k_N(\omega) \right)d\Tilde{R}^\nu(\omega)\leq \frac{1}{(1-\alpha)^{Nd/2}}.
      \end{equation}     
      Finally, for any $\alpha\in (0,1)$, $\nu\in (0,1)$ and $s,t\in [0,1]$ there holds 
            \begin{equation}
      \label{e:k1}
          \int \exp\left(\frac{\alpha\abs{\omega_t-\omega_s}^2}{2\nu\abs{t-s}} \right)d\Tilde{R}^\nu(\omega)\leq \frac{1}{(1-\alpha)^{d/2}}.
      \end{equation}
  \end{lemma}
  
  \begin{proof}
      This lemma is proved in \cite[Lemma 3.3]{baradatm}, where the integration in the left hand side of \eqref{e:knest} is performed with respect to the Reversible Wiener measure on the torus $\mathbb{T}^d$. A quick look at the proof shows that it generalizes to Wiener measures on $\R^d$ with any initial probability distribution (in our case $\Tilde{R}^\nu$, which has initial probability distribution $e^{-\pi\abs{x}^2}dx$). 
      Finally \eqref{e:k1} follows from \eqref{e:knest} and a rescaling.
 \end{proof}
We now provide an existence result. 
\begin{proposition}
\label{prop:existencelagr}
Let $\nu\geq 0$. Suppose \ref{hyp:initial} and either \ref{hyp:mfg} or \ref{hyp:eot} hold. Suppose moreover that there exists $\overline{Q}\in \mathcal{P}(\Omega)$ such that $\mathcal{L}_\nu(\overline{Q})<+\infty$.
     Then the problem 
\eqref{e:problagr}
    admits a minimizer.
\end{proposition}
\begin{proof}
    We prove the result under \ref{hyp:mfg}, since under \ref{hyp:eot} it is classical (see \cite{survey}). We use the direct method of calculus of variations. The existence of $\overline{Q}$ ensures that the infimum is finite. 
    Let us start with the limit problem, that is $\nu=0$. Let $Q^n\in\mathcal{P}(\Omega)$ be a minimizing sequence, so that in particular there exists $C>0$ such that $\mathcal{K}(Q^n)\leq C$ for all $n\in\mathbb{N}$. Consider for $\eps>0$ the set 
    \begin{equation}
        K_\eps=\left\{\omega\in \Omega: \abs{\omega_0}^2+ k(\omega)\leq \frac{1}{\eps} \right\}.
    \end{equation}
    For any $\omega\in AC^{2}([0,1];\R^d)$ and for any $s,t\in [0,1]$ there holds  
    \begin{equation}
    \label{e:equibound}
        \abs{\omega_t}^2\leq 2\abs{\omega_0}^2+4k(\omega)
    \end{equation}
    and
    \begin{equation}
    \label{e:equicon}
        \abs{\omega_t-\omega_s}^2\leq 2k(\omega)\abs{t-s}.
    \end{equation}
    Conditions \eqref{e:equibound} and \eqref{e:equicon} imply that $K_\eps$ is equibounded and equicontinuous, hence by Ascoli-Arzelà's theorem $K_\eps$ is compact. 
    Moreover it holds 
    \begin{equation}
        Q^n(K_\eps^c)\leq \eps\left(\int_{\R^d}\abs{x}^2\,d\rho_0(x)+\int k(\omega)dQ^n(\omega)\right)\leq  \eps\left(\int_{\R^d}\abs{x}^2\,d\rho_0(x)+C\right),
    \end{equation}
    so that $\{Q^n\}_{n\in\N}$ is tight, and Prohorov's Theorem (see for instance \cite{figalli_glaudo}) implies that there exists a (non relabeled) sequence  $\{Q^n\}_{n\in\N}$ converging narrowly to $Q\in \mathcal{P}(\Omega)$. Thus, for any $t\in [0,1]$ there holds $Q^n_t\rightarrow Q_t$ narrowly.
\begin{comment}
     Moreover integrating \eqref{e:equicon} and comparing with the Wasserstein distance we get that for any $s,t\in [0,1]$
    \begin{equation}
        W^2_2(\mu_t,\mu_s)\leq \int \abs{\omega_t-\omega_s}^2dQ_n(\omega)\leq C\abs{t-s},
    \end{equation}
    hence, by Ascoli-Arzelà (on $C^{0}([0,1];\mathcal{P}_2(\R^d)$) we can conclude (again up to a non relabeled subsequence), that $Q^n_t\rightarrow Q_t$ in $C^{0}([0,1];\mathcal{P}_2(\R^d)$.    
\end{comment}
    Existence then follows by the narrow lower semicontinuity of $F,G$ (plus Fatou's lemma) and by the lower semicontinuity of $\mathcal{K}$ with respect to narrow convergence on $\mathcal{P}(\Omega)$. 

    Let us now consider the case $\nu>0$. Consider again a minimizing sequence $\{Q^n\}_{n\in\N}$. By \eqref{e:entrtrick} there holds 
    \begin{equation}
        \int \frac{dQ^n}{d\Tilde{R}^\nu}\log \left(\frac{dQ^n}{d\Tilde{R}^\nu}\right)\,d\Tilde{R}^\nu=H(Q^n|\Tilde{R}^\nu)\leq C+ \pi\int_{\R^d}\abs{x}^2\,d\rho_0(x).
    \end{equation}
    Since $t\rightarrow t \log t$ is superlinear, de la Vallée-Poussin Theorem (see for instance \cite[Theorem 2.29]{fonseca_leoni_calcvar}) implies that $\{\frac{dQ^n}{d\Tilde{R}^\nu}\}_{n\in\N}$ is equi-integrable, and thus, since $\Tilde{R}^\nu$ is a finite measure, that $\{Q^n\}_{n\in\N}$ is tight (see for instance \cite[Lemma 2.1.9]{figalli_glaudo}). Again Prohorov's theorem implies that there exists a non-relabeled subsequence that converges narrowly in $\mathcal{P}(\Omega)$.
\begin{comment}    
    Moreover, by \eqref{e:duality} with $f(\omega)=\frac{\abs{\omega_t-\omega_s}^2}{4\nu\abs{t-s}}$
    \begin{equation}
    \label{e:legagain}
           H(Q_n|\Tilde{R}^{\nu})\geq \int f(\omega) dQ_n(\omega)- \log\int\exp\left(f(\omega)\right)d\Tilde{R}^{\nu}(\omega).
    \end{equation}
    By \eqref{e:k1} with $\alpha=1/2$ the last term in \eqref{e:legagain} is bounded, and we deduce
    \begin{equation}
        \int \frac{\abs{\omega_t-\omega_s}^2}{\abs{t-s}}\, dQ_n(\omega)\lesssim_{d,\nu} 1.
    \end{equation}
\end{comment}
    From here, the proof of existence is analogous to the case $\nu=0$, noticing that $H(\cdot|\Tilde{R}^{\nu})$ is lower semi-continuous with respect to narrow convergence.
\end{proof}
    The next sections are devoted to proving the $\Gamma$-$\liminf$ inequality (Proposition \ref{prop:lower}) and the $\Gamma$-$\limsup$ inequality (Proposition \ref{prop:upper}).
    \subsection{Lower bound and compactness}
  We start by proving an equicompactness result in the same spirit of Proposition \ref{prop:compeul}. For that, we need the following Lemma, which is  a direct consequence of Schilder's theorem (see for instance \cite[Chapter 5]{dembo}).
 \begin{lemma}[Exponential tightness]
 \label{lem:schilder}
     Let $\nu_n\rightarrow 0$. Then for any $L>0$ there exists a compact set $K_L\subset \Omega$ such that
     \begin{equation}
     \label{e:exptight}
         \limsup_n \nu_n \log R^{\nu_n,0}(K_L^c)\leq -L.
     \end{equation}
 \end{lemma}
  \begin{proposition}
  \label{prop:complagr}
      Let $\nu_n\rightarrow 0$, $\rho_0\in \mathcal{P}_2(\R^d)$ with $\ent(\rho_0)<+\infty$ and $\{Q^n\}_{n\in\N}\subset \mathcal{P}(\Omega)$ such that $Q^n_0=\rho_0$ for all $n\in\N$ and $\nu_nH(Q^n|R^{\nu_n})\leq C$. Then there exists a subsequence $Q^{n_k}\rightarrow Q\in \mathcal{P}(\Omega)$ narrowly and such that 
      $(Q^{n_k}_t)_{t\in [0,1]}\rightarrow (Q_t)_{t\in [0,1]}\in C^{0}([0,1],\mathcal{P}_2(\R^d))$ narrowly.
  \end{proposition}
  \begin{proof}
      By \eqref{e:entrtrick}, we also have \begin{equation}
          \label{e:nuentrnu}\nu_nH(Q^n|\Tilde{R}^{\nu_n})\leq C.
      \end{equation} We apply Proposition \ref{prop:entrdisint} with respect to $e_0$ to obtain 
        \begin{equation}
        \label{e:disin}H(Q^n|\Tilde{R}^{\nu_n})=H(\rho_0|e^{-\pi\abs{x}^2}dx)+\int_{\R^d}H(Q^{n,x}|R^{\nu_n,x})\,d\rho_0(x).
        \end{equation}
        Define $\hat{Q}^{n,x}=(\cdot-x)\#Q^{n,x}$, so that 
        \begin{equation}
        \label{e:shift}H(Q^{n,x}|R^{\nu_n,x})=H(\hat{Q}^{n,x}|R^{\nu_n,0}).
        \end{equation}
      Let $L>0$ and choose $K_{2L}$ as given by Lemma \ref{lem:schilder}. Define the function $f(\omega)=\frac{L}{\nu_n}\chi_{K_{2L}^c}(\omega)$. We test \eqref{e:duality} with $f$ and multiply by $\nu_n$, obtaining 
      \begin{equation}  \nu_nH(\hat{Q}^{n,x}|R^{\nu_n,0})\geq L \hat{Q}^{n,x}(K_{2L}^c)-\nu_n\log(1+(e^{\frac{L}{\nu_n}}-1)R^{\nu_n,0}(K_{2L}^c)). 
      \end{equation}
      By \eqref{e:exptight}, $R^{\nu_n,0}(K_{2L}^c)\leq e^{-3L/(2\nu_n)}$ for $n$ large, so the last term tends to $0$, and we can conclude that 
      \begin{equation}
      \label{e:xexptight}
    \hat{Q}^{n,x}(K_{2L}^c)\leq \frac{\nu_n}{L}H(\hat{Q}^{n,x}|R^{\nu_n,0})+\frac{C}{L}.
      \end{equation}
      We now integrate \eqref{e:xexptight} and find using \eqref{e:disin}, \eqref{e:shift} and \eqref{e:nuentrnu}
      \begin{equation}
      \label{e:partialtight}
          Q^n(\{\omega-\omega_0\in K_{2L}^c\})= \int_{\R^d}\hat{Q}^{n,x}(K_{2L}^c)d\rho_0(x)\leq \frac{\nu_n}{L}\int_{\R^d}H(\hat{Q}^{n,x}|R^{\nu_n,0})d\rho_0(x)+\frac{C}{L}\leq \frac{C}{L}.
      \end{equation} 
      We finally define \begin{equation}
          \Tilde{K}_L=\{\omega
          \in \Omega: \abs{\omega_0}\leq L,\, \omega-\omega_0 \in K_{2L}\}.
    \end{equation} 
    $\Tilde{K}_L$ is compact since $K_{2L}$ is compact. Using \eqref{e:partialtight} we obtain
    \begin{multline}
    Q^n(\Tilde{K}_L^c)\leq Q^n(\{\abs{\omega_0}>L\})+ Q^n(\{\omega-\omega_0\in K_{2L}^c\})\\ = \rho_0(\{\abs{x}>L\})+ Q^n(\{\omega-\omega_0\in K_{2L}^c\})\leq \frac{C
    }{L},
     \end{multline}
     where in the last passage we also use the finiteness of the first moment of $\rho_0$. We have thus proved tightness, hence
    Prohorov's theorem implies the existence of a non-relabeled subsequence $Q^n\rightarrow Q\in \mathcal{P}(\Omega)$ narrowly. Convergence of the marginals is given by Lemma \ref{lem:cont2}, we are left to prove that $(Q_t)_{t\in [0,1]}\in C^{0}([0,1],\mathcal{P}_2(\R^d))$. Let $t,s\in[0,1]$. Again by \eqref{e:duality}, this time with $f(\omega)=\frac{\abs{\omega_t-\omega_s}^2}{4\nu_n\abs{t-s}}$
    \begin{equation}
    \label{e:legagain}
           C\geq \nu_n H(Q_n|\Tilde{R}^{\nu_n})\geq \int \nu_n f(\omega) dQ_n(\omega)- \nu_n\log\int\exp\left(f(\omega)\right)d\Tilde{R}^{\nu_n}(\omega).
    \end{equation}
    By \eqref{e:k1} with $\alpha=1/2$ the last term in \eqref{e:legagain} is bounded, and we deduce
    \begin{equation}
        \int \frac{\abs{\omega_t-\omega_s}^2}{\abs{t-s}}\, dQ_n(\omega)\lesssim_{d} 1,
    \end{equation} 
  hence by using the competitor $\pi_n=(e_t,e_s)\#Q_n$
  \begin{equation}
  \label{e:equiconlagr}
      W_2^2(Q_t^n,Q_s^n)\leq \int_{\R^d\times\R^d}\abs{x-y}^2d\pi_n\leq \int \abs{\omega_t-\omega_s}^2\, dQ_n(\omega)\lesssim_d\abs{t-s}.
  \end{equation}
  Using \eqref{e:equiconlagr}, the same argument as in Proposition \ref{prop:compeul} finishes the proof.
\begin{comment}
  Finally, 
  \begin{equation}
      \int_{\R^d}\abs{x}^2dQ_t^n=W_2^2(Q_t^n,\delta_0)\lesssim W_2^2(\rho_0,\delta_0)+W_2^2(\rho_0,Q_t^n)\overset{\eqref{e:equicon}}{\lesssim_d}1,  
  \end{equation}
  hence Ascoli-Arzelà theorem on $C^0([0,1],\mathcal{P}(\R^d))$ with the metric \eqref{e:narrow} (see \cite[Proposition 3.3.1]{AGS}) gives the convergence of the marginals. Finally we are left to prove that $(Q_t)_{t\in[0,1]}\in C^0([0,1],\mathcal{P}_2(\R^d))$, but this is a direct consequence of \eqref{e:equicon} since $W_2^2(\cdot,\cdot)$ is lower semicontinuous with respect to narrow converngence.
\end{comment}
  \end{proof}
  We now prove the $\Gamma$-$\liminf$ inequality.   \begin{proposition}
  \label{prop:lower}
      Suppose \ref{hyp:initial} and either \ref{hyp:mfg} or \ref{hyp:eot} hold. Let $\nu_n\rightarrow 0$. Let $Q_n,Q\in \mathcal{P}(\Omega)$ such that $Q_n\rightarrow Q$ narrowly. Then
       \begin{equation}
           \liminf_n\mathcal{L}_{\nu_n}(Q_n)\geq \mathcal{L}_0(Q).
       \end{equation}
  \end{proposition}
  This result could be derived from \cite{leonard_convergence}, and is a minor adaptation of \cite[Theorem 3.1]{baradatm}. For completeness, we provide the argument. 
  \begin{proof}[Proof of Proposition \ref{prop:lower}]
      As mentioned, the proof is a minor adaptation of \cite[Theorem 3.1]{baradatm}. We provide the main steps of the proof. 
      
      Let us remark first that since by Lemma \ref{lem:cont2} there holds  $(Q^{n}_t)_{t\in[0,1]}\rightarrow (Q_t)_{t\in[0,1]}$ narrowly in $ C([0,1];\mathcal{P}(\R^d))$, we can follow for $C$ the same strategy as in Theorem \ref{th:maineul} (under \ref{hyp:eot}, $C\equiv 0$ for configurations with finite energy). 
      
      For the entropy term, recall that by \eqref{e:entrtrick} there holds $H(Q^n|R^\nu)=H(Q^n|\Tilde{R}^\nu)-\pi\int_{\R^d}\abs{x}^2\,d\rho_0(x)$. We use \eqref{e:duality} with $f(\omega)=\alpha k_N(\omega)/\nu_n$ where $\alpha\in (0,1)  $ and $k_N$ is defined in Definition \ref{def:kn}. There holds 
      \begin{equation}
      \label{e:legendreest}
          \nu_n H(Q_n|\Tilde{R}^{\nu_n})\geq \int \alpha k_N(\omega) dQ_n(\omega)-\nu_n\log \int\exp\left(\frac{\alpha}{\nu_n}k_N(\omega)\right)d\Tilde{R}^{\nu_n}(\omega)
      \end{equation}
By \eqref{e:knest}, the last term in \eqref{e:legendreest} is a bounded quantity multiplied by $\nu_n$, hence
\begin{equation}
    \liminf_n \nu_n H(Q_n|R^{\nu_n})=\liminf_n \nu_n H(Q_n|\Tilde{R}^{\nu_n}) \geq \liminf_n \int \alpha k_N(\omega) dQ_n(\omega)= \alpha \int  k_N(\omega)dQ(\omega),
\end{equation}
      where in the last passage we used narrow convergence of $Q_n$ paired with continuity and boundedness of $k_N$. The conclusion then follows by sending $\alpha\rightarrow 1$ and $N\rightarrow \infty$, employing  Fatou's lemma and the pointwise convergence of $k_N$.
  \end{proof}
  \subsection{Upper bound}
  The $\Gamma$-$\limsup$ inequality is more involved. 
\begin{comment}
      First we need a density result.
  \begin{proposition}
      Let $Q\in P(\Omega)$ with $\mathcal{F}_0(Q)<+\infty$. There exist $Q_n\rightarrow Q$ narrowly in $\Omega$ such that $Q_n$ is "piecewise constant" and  
      \begin{equation}
          \limsup_n\mathcal{F}_0(Q_n)\leq \mathcal{F}_0(Q).
      \end{equation}
  \end{proposition}
  \begin{proof}
      Idea: piecewise approximation. I think we fix $N$ equispaced marginals, interpolate (how?) to find a flow of measures and then implicitly define $Q^N$ using \cite[Theorem 8.3.1]{AGS}?
  \end{proof}
  Once we have the density result, we can restrict the construction of recovery sequences to paths with finitely many prescribed marginals (multimarginal problem).
\end{comment}
  To understand the strategy, we first prove the result for \ref{hyp:eot}, that is an alternative proof of the recovery sequence for the dynamic entropic optimal transport problem. An appropriate adaptation then will work for the setting of \ref{hyp:mfg}.
  In the following, for simplicity, we denote for $\gamma\in \mathcal{P}(\R^d\times\R^d)$
\begin{equation}
    c(\gamma)=\frac{1}{2}\int_{\R^d\times \R^d}\abs{x-y}^2d\gamma(x,y).
\end{equation}
  This first proposition is a simpler case of \cite[Theorem 3.1]{baradatm}. 
  \begin{proposition}
  \label{prop:reconstruct}
   Let $Q\in \mathcal P(\Omega)$ such that $\mathcal{L}_0(Q)<+\infty$, $\ent(\rho_0)<+\infty$ and $\ent(\rho_1)<+\infty$. Consider $\gamma=(e_0,e_1)\#Q$.
      Let $\nu_n\rightarrow 0$ and $\gamma_n\in \mathcal{P}(\R^d\times\R^d)$ such that $W_2(\gamma_n, \gamma)\rightarrow 0$ and 
      \begin{equation}
      \label{e:staticlimsup}
          \limsup_n\nu_n H(\gamma_n|R_{0,1}^{\nu_n})\leq c(\gamma).
      \end{equation}
      Then, there exists $(Q_n)_{n\in \N}\subset \mathcal P(\Omega)$ such that $W_2(Q_n,Q)\rightarrow 0$, $\gamma_n=(e_0,e_1)\#Q_n$ and
      \begin{equation}
          \limsup_n{\nu_nH(Q_n|R^{\nu_n})}\leq \mathcal{K}(Q).
      \end{equation}
  \end{proposition}
To prove this we need some preliminary results. 
    Let $B^{\nu}=R^{\nu,0,0}$ be the Brownian Bridge connecting $0$ to $0$ with diffusivity $\nu$. For $\omega\in AC^2([0,1],\R^d)$, define the map $T_\omega: \Omega\rightarrow  \Omega$ by
\begin{equation}
    T_\omega( \alpha)=\alpha+\omega
\end{equation}
and the Brownian Bridge "around" $\omega$ as 
\begin{equation}
B_\omega^{\nu}=T_{\omega}\#B^{\nu}.
\end{equation}
\begin{lemma}
\label{lem:CM}
    For all $\nu>0$ and $\omega\in AC^2([0,1],\R^d)$
    \begin{equation}
        \nu H(B_{\omega}^{\nu}|R^{\nu,\omega_0,\omega_1})=k(\omega)-\frac{1}{2}\abs{\omega_1-\omega_0}^2.
    \end{equation}
\end{lemma}
\begin{proof}
    See \cite[Lemma B.2]{baradatm}.
\end{proof}
\begin{proposition}
\label{prop:entrbridge}
    Let $P\in \mathcal{P}(\Omega)$ such that $\mathcal{K}(P)<+\infty$ and $0<\nu<1$. Define 
    \begin{equation}
        Q=\int B_{\omega}^{\nu}dP(\omega).
    \end{equation}
    Then $(e_0,e_1)\#Q=(e_0,e_1)\#P=\gamma$ and
    \begin{equation}
\label{e:margconv}Q_t=P_t\ast\tau_{\nu t(1-t)}.
    \end{equation}
    Moreover 
    \begin{equation}
    \label{e:wassbridge}
        W_2^2(P,Q)\lesssim_d \nu
    \end{equation}
    and 
    \begin{equation}
    \label{e:entrkin}
\nu H(Q|R^{\nu})\leq \nu H(\gamma|R^{\nu}_{0,1})+\mathcal{K}(P)-c(\gamma).
    \end{equation}
\end{proposition}
\begin{proof}
    Q can be thought of as a "noisy" version of $P$. The endpoints are left unchanged, hence $(e_0,e_1)\#Q=(e_0,e_1)\#P$. For \eqref{e:margconv} and\eqref{e:wassbridge} we use a probabilistic argument. Consider two independent random variables $X\sim P$ and $\beta^\nu\sim B^\nu$. By definition \begin{equation}
        Y=X+\beta^\nu\sim Q.
    \end{equation}
    By \eqref{e:bridge}, \eqref{e:margconv} holds. Moreover, if we consider the coupling $(X,Y)$
    \begin{equation}
        W_2^2(P,Q)\leq \mathbb{E}\left[\sup_t\abs{Y_t-X_t}^2\right]=\mathbb{E}\left[\sup_t\abs{\beta_t^\nu}^2\right]\lesssim_d \nu,
    \end{equation}
    where in the last passage we use \eqref{e:linfest}.
    Finally we prove \eqref{e:entrkin}. By Proposition \ref{prop:entrdisint}
    \begin{equation}
        \label{e:cond}H(Q|R^{\nu})=H(\gamma|R^{\nu}_{0,1})+\int H(Q^{x,y}|R^{\nu,x,y})\,d\gamma(x,y)
    \end{equation}
    and 
    \begin{equation}
        Q^{x,y}=\int B^{\nu}_\omega \,dP^{x,y}(\omega).
    \end{equation}
    By convexity of the relative entropy (Jensen's inequality) and remarking that $P^{x,y}$-almost surely $\omega_0=x$ and $\omega_1=y$
    \begin{equation}
    \label{e:jensen}
        H(Q^{x,y}|R^{\nu,x,y})\leq \int H(B^{\nu}_\omega|R^{\nu,x,y})\,dP^{x,y}(\omega)=\int H(B^{\nu}_\omega|R^{\nu,\omega_0,\omega_1})\,dP^{x,y}(\omega).
    \end{equation}
    Hence, combining \eqref{e:cond} and \eqref{e:jensen} we get 
    \begin{multline}
        H(Q|R^{\nu})\leq H(\gamma|R^{\nu}_{0,1})+\int \int H(B^\nu_\omega|R^{\nu,\omega_0,\omega_1})dP^{x,y}(\omega)d\gamma(x,y)\\ =H(\gamma|R^{\nu}_{0,1})+\int H(B^\nu_\omega|R^{\nu,\omega_0,\omega_1})dP(\omega).
     \end{multline}
      Multiplying by $\nu$ and exploiting Lemma \ref{lem:CM}, we obtain the result.
\end{proof}
We now prove a stability result for the kinetic energy.
\begin{proposition}
\label{prop:stability}
    Let $P\in \mathcal{P}(\Omega)$, $\gamma=(e_0,e_1)\#P\in \mathcal{P}(\R^d\times \R^d)$ and $\eta\in\mathcal{P}(\R^d\times \R^d)$. There exists $Q\in \mathcal{P}(\Omega)$ with $\eta=(e_0,e_1)\#Q$, 
    \begin{equation}
        W_2^2(P,Q)\lesssim W_2^2(\gamma,\eta)
    \end{equation}
    and 
    \begin{equation}
        \mathcal{K}(Q)\leq \mathcal{K}(P)+2\sqrt{\mathcal{K}(P)}W_2(\gamma,\eta)+ W_2^2(\gamma,\eta).
    \end{equation}
\end{proposition}
\begin{proof}
    The proof is an adaptation of \cite[Lemma 5]{baradatcont}. Let $X\sim P$, so that $(X_0,X_1)\sim \gamma$. Consider $(Y_0,Y_1)\sim \eta$ such that $((X_0,X_1),(Y_0,Y_1))$ is optimal for \eqref{e:wassprob}. We define 
    \begin{equation}
        Y_t=X_t +(1-t)(Y_0-X_0)+t(Y_1-X_1)
    \end{equation}
     and $Q=\operatorname{Law}(Y)$. By construction $\eta=(e_0,e_1)\#Q$.
     Notice that 
     \begin{equation}
         d_{\Omega}(X,Y)=\sup_{t\in [0,1]}\abs{X_t-Y_t}\leq \abs{X_0-Y_0}+\abs{X_1-Y_1},
     \end{equation}
     hence 
     \begin{equation}
         W_2^2(P,Q)\leq \mathbb{E}[d_\Omega(X,Y)^2]\leq \mathbb{E}[(\abs{X_0-Y_0}+\abs{X_1-Y_1})^2]\leq 2 W_2^2(\gamma,\eta).
     \end{equation}
     Finally let us compute the kinetic energy.   Recall that if $X$ is an $\Omega$-valued random variable with law $Q$, then $\mathcal{K}(Q)= \mathbb E[k(X)]$.
One has
\begin{multline}
    k(Y)=\frac{1}{2}\int_{0}^1\abs{\dot{Y}(t)}^2\,dt=\frac{1}{2}\int_{0}^1\abs{\dot{X}(t)-(Y_0-X_0)+(Y_1-X_1)}^2\,dt  
    \\ \leq \frac{1}{2}\left(\norm{\dot{X}}_{L^2(0,1)}+\abs{Y_0-X_0}+\abs{Y_1-X_1}\right)^2 \leq \frac{1}{2}\left(\norm{\dot{X}}_{L^2(0,1)}+\sqrt{2}\abs{(Y_0,Y_1)-(X_0,X_1)}\right)^2\\=\left( \sqrt{k(X)}+\abs{(Y_0,Y_1)-(X_0,X_1)}\right)^2,
\end{multline}
so that,
\begin{multline}
        \mathcal{K}(Q)=\mathbb{E}[k(Y)] \leq \mathbb{E}\left[\left( \sqrt{k(X)}+\abs{(Y_0,Y_1)-(X_0,X_1)}\right)^2\right]\\ \leq
        \left(\sqrt{\mathbb{E}\left[k(X)\right]}+\sqrt{\mathbb{E}\left[\abs{(Y_0,Y_1)-(X_0,X_1)}^2\right]}\right)^2\\
        = \left(\sqrt{K(P)}+W_2(\gamma,\eta)\right)^2=\mathcal{K}(P)+2\sqrt{\mathcal{K}(P)}W_2(\gamma,\eta)+W_2^2(\gamma,\eta).
    \end{multline}
\end{proof}
\begin{proof}[Proof of Proposition \ref{prop:reconstruct}]
   This proof is an easier case of \cite[Theorem 3.1]{baradatm}. Consider $\gamma_n\rightarrow \gamma$. By applying Proposition \ref{prop:stability}, there exists $(P_n)_{n\in \N}\subset P(\Omega)$ such that $W_2(P_n,Q)\rightarrow 0$, $(e_0,e_1)\#P_n=\gamma_n$ and $\limsup_n\mathcal{K}(P_n)\leq \mathcal{K}(Q)$. We define 
    \begin{equation}
        Q_n=\int B_{\omega}^{\nu_n}dP^n(\omega).
    \end{equation}
    By \eqref{e:wassbridge} it follows $W_2(Q_n,Q)\rightarrow 0$ and by \eqref{e:margconv} the marginal constraints hold. Applying \eqref{e:entrkin}, 
    \begin{equation}
    \label{e:est}
\nu_n H(Q_n|R^{\nu_n})\leq \nu_n H(\gamma_n|R^{\nu_n}_{0,1})+\mathcal{K}(P_n)-c(\gamma_n).
    \end{equation}
     Taking the limsup in \eqref{e:est}, the result follows.
\end{proof}
We now adapt Proposition \ref{prop:reconstruct} to the setting \ref{hyp:mfg}. Here, Proposition \ref{prop:stability} is replaced by the following ad-hoc stability result.
\begin{proposition}
\label{prop:stability2}
        Let $P\in \mathcal{P}(\Omega)$, $\gamma=(e_0,e_1)\#P\in \mathcal{P}(\R^d\times \R^d)$, with $P_0=\rho_0$. Consider the disintegration $\gamma(dx,dy)=\gamma_x(dy)\otimes\rho_0(dx)$. Consider for $\rho_0-$a.e.\ $x$ the measure $\gamma^\nu_x=\gamma_x\ast \tau_{\nu}$. Define $\gamma^\nu(dx,dy)=\gamma^\nu_x(dy)\otimes\rho_0(dx)$. There exists  $P^\nu \in \mathcal{P}(\Omega)$ with $\gamma^\nu=(e_0,e_1)\#P^\nu$ and such that for all $t\in [0,1]$
        \begin{equation}
        \label{e:margprop}
            P^\nu_t=P_t\ast \tau_{\nu t^2}.
        \end{equation}
        Moreover 
    \begin{equation}
    \label{e:wasspaths}W_2^2(P,P_\nu)\leq d\nu
    \end{equation}
    and 
    \begin{equation}
    \label{e:kinpaths}
        \mathcal{K}(P_\nu)= \mathcal{K}(P)+\frac{d\nu}{2}.
    \end{equation}
\end{proposition}
\begin{proof}
    Take a random variable $X\sim P$ and an independent variable $Z_\nu\sim\mathcal{N}(0,\nu I_d)$. Define for $t\in [0,1]$ the random variable
    \begin{equation}
        Y_t=X_t+tZ_\nu
    \end{equation}
and $P_\nu=\operatorname{Law}(Y)$. \eqref{e:margprop} follows by definition. For \eqref{e:wasspaths}, consider the coupling $(X,Y)$. It holds 
\begin{equation}
    W_2^2(P,P_\nu)\leq \mathbb{E}\left[\sup_t\abs{Y_t-X_t}^2\right]=\mathbb{E}[\abs{Z_\nu}^2]=d\nu.
\end{equation}
For \eqref{e:kinpaths}, one has using independence and centering of $Z$
\begin{equation}
\mathcal{K}(P_\nu)=\frac12\mathbb{E}\left[\int_0^1\abs{\dot{X}_t+Z_\nu}^2dt\right]=\mathcal{K}(P)+\frac12\mathbb{E}\left[\abs{Z_\nu}^2\right]=\mathcal{K}(P)+\frac{d\nu}{2}.
\end{equation}
\end{proof}

 Recall that, by a change of reference measure, for any $\gamma\in \mathcal{P}(\R^d\times\R^d)$ for which the below terms are defined   
  \begin{equation}
  \label{e:varchange}
      \nu H(\gamma|R_{0,1}^{\nu})=\nu \ent(\gamma)+c(\gamma)+\frac{d}{2}\nu\log(2\pi\nu).
  \end{equation}
  \begin{comment}
   We introduce the so called block approximation of entropic optimal transport from \cite{carlier2017convergence}.
  The following is taken by \cite[Proposition 2.14]{carlier2017convergence}.
\begin{lemma}
\label{lem:block}
    Let $\rho_0,\rho_1\in \mathcal{P}_2(\R^d)$, $H(\rho_0)<+\infty$,  $H(\rho_1)<+\infty$, $\gamma\in \Pi(\rho_0,\rho_1)$ and $0<l<1$. There exists $\gamma_l\in \Pi(\rho_0,\rho_1)$ and numbers $\alpha=\alpha(d)\in (0,1),\, C=C(d)>0$ such that 
    \begin{equation}
    \label{e:blockest}
        \ent(\gamma_l)\leq \ent(\rho_0)+\ent(\rho_1)+C\left((M(\rho_0)+dl^2+1)^\alpha+(M(\rho_1)+d\nu_n^2+1)^\alpha-2d\log l\right)
    \end{equation}
    and 
    \begin{equation}
        W_2^2(\gamma,\gamma_l)\leq 2dl^2.
    \end{equation}
\end{lemma}
\end{comment}
We finally prove the $\Gamma$-$\limsup$ inequality. 
\begin{proposition}
\label{prop:upper}
    Suppose \ref{hyp:initial} and either \ref{hyp:mfg} or \ref{hyp:eot} hold. Let $\nu_n\rightarrow 0$. Let $Q\in \mathcal{P}(\Omega)$ such that $\mathcal{L}_0(Q)<+\infty$. There exists $(Q^n)_{n\in \N}\subset \mathcal P(\Omega)$ such that $W_2(Q^n,Q)\rightarrow 0$ and
       \begin{equation}
           \limsup_n\mathcal{L}_{\nu_n}(Q^n)\leq \mathcal{L}_0(Q).
       \end{equation}
\end{proposition}
\begin{proof}
    Let us start by the case of \ref{hyp:eot}. Consider $\gamma=(e_0,e_1)\#Q$. By \cite{carlier2017convergence}, there exists a sequence $\{\gamma_n\}_{n\in \N}\subset \Pi(\rho_0,\rho_1)$ with $W_2(\gamma_n,\gamma)\to 0$ such that \eqref{e:staticlimsup} holds. Proposition \ref{prop:reconstruct} then gives the result. 
    
    We now consider \ref{hyp:mfg}. Consider $\gamma=(e_0,e_1)\#Q$. By disintegration we can write $\gamma(dx,dy)=\gamma^x(dy)\otimes\rho_0(dx)$. Consider for $\rho_0-$a.e.\ $x$ the measure $\gamma^x_n=\gamma^x\ast \tau_{\nu_n}$. Define $\gamma_n(dx,dy)=\gamma^x_n(dy)\otimes\rho_0(dx)$. By Lemma \ref{lem:moll}, $e_1\#\gamma_n=Q_1\ast\tau_{\nu_n}$ and $W_2(\gamma_n,\gamma)\rightarrow 0$, so that $c(\gamma_n)\rightarrow c(\gamma)$ by classical properties of optimal transport, see for instance \cite[Theorem 6.9]{Villani}. Now we can repeat the same argument as in Proposition \ref{prop:reconstruct}.
    By Proposition \ref{prop:stability2}, there exists $(P^n)_{n\in \N}\subset P(\Omega)$ such that $W_2(P^n,Q)\rightarrow 0$, $(e_0,e_1)\#P^n=\gamma_n$ and $\mathcal{K}(P^n)\rightarrow \mathcal{K}(Q)$. We define 
    \begin{equation}
        Q^n=\int B_{\omega}^{\nu_n}dP^n(\omega),
    \end{equation}
   By \eqref{e:wassbridge} and the triangle inequality it follows $W_2(Q^n,Q)\rightarrow 0$, while by \eqref{e:margconv} the marginal constraints hold. Applying \eqref{e:entrkin}, 
    \begin{equation}
    \label{e:est2}
\nu_n H(Q^n|R^{\nu_n})\leq \nu_n H(\gamma_n|R^{\nu_n}_{0,1})+\mathcal{K}(P^n)-c(\gamma_n).
    \end{equation}
    Combining with \eqref{e:varchange}
    \begin{multline}
    \label{e:entrestgamma}
         \nu_n H(Q^n|R^{\nu_n})\leq \nu_n H(\gamma_n|R^{\nu_n}_{0,1})+\mathcal{K}(P^n)-c(\gamma_n)\\ \overset{\eqref{e:varchange}}{=} \mathcal{K}(P^n)+ \nu_n \ent(\gamma_n)+\frac{d}{2}\nu_n\log(2\pi\nu_n).
    \end{multline}
  We now estimate $\ent(\gamma_n)$. By \eqref{e:entrgamma}
  \begin{equation}
      \ent(\gamma_n)\leq  \ent(\rho_0)-\frac{d}{2}\log(2\pi \nu_n).
  \end{equation}
  This implies that $\limsup_n\nu_n\ent(\gamma_n)\leq 0$. By taking the $\limsup$ in \eqref{e:entrestgamma}, we obtain 
  \begin{equation}
      \limsup \nu_n H(Q^n|R^{\nu_n})\leq\mathcal{K}(Q).
  \end{equation}
To conclude, let us prove continuity of $C$. By \eqref{e:margconv}, \eqref{e:margprop} and properties of the convolution, it holds 
\begin{equation}
    Q^n_t=Q_t\ast \tau_{\nu_n t^2}\ast\tau_{\nu_n t(1-t)}=Q_t\ast\tau_{\nu_n t}.
\end{equation}
The proof of the continuity of $C$ can thus be replicated as in Theorem \ref{th:maineul} (note that $(Q^n_t)_t$ is exactly the Eulerian recovery sequence used there).
\end{proof}

\section*{Acknowledgements and AI Statement}
A.C.\ and L.N. are partially supported by the ANR project
GOTA (ANR-23-CE46-0001). The authors acknowledge the use of AI tools. In particular, ChatGPT 6 was used to check proofs, help with the writing and ensure consistency of notation. ChatGPT 6 also helped with the proofs of Theorem \ref{th:energyexp} and Proposition \ref{prop:complagr}, and produced the examples in Section \ref{sec:eul}. The authors have verified all results and the text, and take full responsibility for the final manuscript.
\printbibliography
\end{document}